\documentclass{article}
 
\usepackage{amsmath,amssymb,amsfonts}
\usepackage{pstricks, pst-node}
\usepackage[all]{xy}
\usepackage[dvips]{graphics}

\def\Id{\mathrm{Id}}
\def\id{\mathrm{id}}
\def\Hom{\mathrm{Hom}}

\def\N{\mathbb{N}}

\def\Q{\mathbb{Q}}

\def\M{{\cal{M}}}
\def\MM{\overline{\cal{M}}}

\def\P{{\cal{P}}}
\def\PP{\overline{\cal{P}}}
\def\QQ{{\cal{Q}}}
\def\F{{\cal{F}}}
\def\Endom{{\cal{E}}nd}
\def\As{{\cal{A}}s}
\def\Dias{{\cal{D}}ias}

\def\Dup{{\cal{D}}up}

\def\Y{\overline{Y}}
\def\YY{\overline{\overline{Y}}}

\def\over{\slash}
\def\under{\backslash}

\def\io{i^{\mathrm{over}}}
\def\iu{i^{\mathrm{under}}}

\def\Gi{G^{\mathrm{inv}}}
\def\Gp{G^{\gamma}}
\def\Ge{G^{e}}
\def\Go{G^{\mathrm{over}}}
\def\Gu{G^{\mathrm{under}}}
\def\Gd{G^{\mathrm{dif}}}
\def\Gl{G^{\mathrm{\lambda}}}
\def\Gr{G^{\mathrm{\rho}}}
\def\Ga{G^{\alpha}}

\def\H{{\cal{H}}}
\def\Hnc{{\cal{H}}^{\mathrm{nc}}}
\def\Hd{{\cal{H}}^{\mathrm{dif}}}
\def\Hdnc{{\cal{H}}^{\mathrm{dif,nc}}}
\def\Hi{{\cal{H}}^{\mathrm{inv}}}
\def\Hinc{{\cal{H}}^{\mathrm{inv,nc}}}
\def\Hl{{\cal{H}}^{\mathrm{\lambda}}}

\def\Hr{{\cal{H}}^{\mathrm{\rho}}}
\def\Hrnc{{\cal{H}}^{\mathrm{\rho,nc}}}
\def\Hp{{\cal{H}}^{\gamma}}
\def\He{{\cal{H}}^{e}}
\def\Ha{{\cal{H}}^{\alpha}}
\def\Hanc{{\cal{H}}^{\alpha,nc}}

\def\Dd{\Delta^{\mathrm{dif}}}
\def\Di{\Delta^{\mathrm{inv}}}
\def\Dl{\Delta^{\ltimes}}
\def\Dip{\Delta^{\mathrm{inv}}_{\gamma}}
\def\Die{\Delta^{\mathrm{inv}}_e}
\def\Da{\Delta^{\alpha}}

\def\Dr{\Delta^{\rho}}

\def\dinv{\delta^{\mathrm{inv}}}
\def\ddif{\delta^{\mathrm{dif}}}
\def\da{\delta^{\alpha}}

\def\dr{\delta^{\rho}}

\newtheorem{theorem}{Theorem}[section]
\newtheorem{proposition}[theorem]{Proposition}
\newtheorem{corollary}[theorem]{Corollary}
\newtheorem{lemma}[theorem]{Lemma}
\newtheorem{defin}[theorem]{Definition}

\newtheorem{rem}[theorem]{Remark}
\newenvironment{remark}{\begin{rem} \em}{\end{rem}}
\newtheorem{exa}[theorem]{Example}

\newenvironment{proof}{\noindent{\em Proof.\/}}{\hfill$\Box$\par\vspace{.2cm}}
\newenvironment{proof of}[2]{\bigskip \noindent{\em Proof of #1~\ref{#2}.\/}}
        {\hfill$\square$\par\vspace{.2cm}}
\numberwithin{equation}{section}
\numberwithin{figure}{section}

\newcommand{\treeO}{
\setlength{\unitlength}{3pt}
\psset{unit=3pt}
\psset{runit=2pt}
\psset{linewidth=0.2}
\begin{pspicture}(0,-.5)(2.5,2)
\psline(1,-1)(1,1.5)
\end{pspicture}}

\newcommand{\treeA}{
\setlength{\unitlength}{3pt}
\psset{unit=3pt}
\psset{runit=2pt}
\psset{linewidth=0.2}
\begin{pspicture}(0,-.5)(2.5,2)
\psline(1,-1)(1,.5)
\psline(1,.5)(0,1.5)
\psline(1,.5)(2,1.5)
\end{pspicture}}

\newcommand{\treeAsmall}{ 
\setlength{\unitlength}{2pt}
\psset{unit=2pt}
\psset{runit=1.2pt}
\psset{linewidth=0.2}
\begin{pspicture}(0,-.5)(2.5,2)
\psline(1,-1)(1,.5)
\psline(1,.5)(0,1.5)
\psline(1,.5)(2,1.5)
\end{pspicture}}

\newcommand\treeAB{
\setlength{\unitlength}{3pt}
\psset{unit=3pt}
\psset{runit=2pt}
\psset{linewidth=0.2}
\begin{pspicture}(0,0)(5,3)
\psline(3,-1)(3,.5)
\psline(3,.5)(1,2.5)
\psline(3,.5)(4,1.5)
\psline(2,1.5)(3,2.5)
\end{pspicture}}

\newcommand\treeBA{
\setlength{\unitlength}{3pt}
\psset{unit=3pt}
\psset{runit=2pt}
\psset{linewidth=0.2}
\begin{pspicture}(0,-.5)(5,3)
\psline(2,-1)(2,.5)
\psline(2,.5)(4,2.5)
\psline(2,.5)(1,1.5)
\psline(3,1.5)(2,2.5)
\end{pspicture}}

\newcommand\treeABC{
\setlength{\unitlength}{3pt}
\psset{unit=3pt}
\psset{runit=2pt}
\psset{linewidth=0.2}
\begin{pspicture}(0,-.5)(5,4.5)
\psline(3,-1)(3,.5)
\psline(3,.5)(0,3.5)
\psline(3,.5)(4,1.5)
\psline(2,1.5)(3,2.5)
\psline(1,2.5)(2,3.5)
\end{pspicture}}

\newcommand\treeBAC{
\setlength{\unitlength}{3pt}
\psset{unit=3pt}
\psset{runit=2pt}
\psset{linewidth=0.2}
\begin{pspicture}(0,-.5)(5,4.5)
\psline(3,-1)(3,.5)
\psline(3,.5)(1,2.5)
\psline(3,.5)(4,1.5)
\psline(2,1.5)(4,3.5)
\psline(3,2.5)(2,3.5)
\end{pspicture}}

\newcommand\treeACA{
\setlength{\unitlength}{3pt}
\psset{unit=3pt}
\psset{runit=2pt}
\psset{linewidth=0.2}
\begin{pspicture}(0,-.5)(6,3.5)
\psline(3,-1)(3,.5)
\psline(3,.5)(0.5,3)
\psline(3,.5)(5.5,3)
\psline(1.5,2)(2.5,3)
\psline(4.5,2)(3.5,3)
\end{pspicture}}

\newcommand\treeCAB{
\setlength{\unitlength}{3pt}
\psset{unit=3pt}
\psset{runit=2pt}
\psset{linewidth=0.2}
\begin{pspicture}(0,-.5)(5,4.5)
\psline(2,-1)(2,.5)
\psline(2,.5)(4,2.5)
\psline(2,.5)(1,1.5)
\psline(3,1.5)(1,3.5)
\psline(2,2.5)(3,3.5)
\end{pspicture}}

\newcommand\treeCBA{
\setlength{\unitlength}{3pt}
\psset{unit=3pt}
\psset{runit=2pt}
\psset{linewidth=0.2}
\begin{pspicture}(0,-.5)(5,4.5)
\psline(2,-1)(2,.5)
\psline(2,.5)(5,3.5)
\psline(2,.5)(1,1.5)
\psline(3,1.5)(2,2.5)
\psline(4,2.5)(3,3.5)
\end{pspicture}}

\newcommand\treeABCD{
\setlength{\unitlength}{3pt}
\psset{unit=3pt}
\psset{runit=2pt}
\psset{linewidth=0.2}
\begin{pspicture}(0,-.5)(6,5)
\psline(4,-1)(4,.5)
\psline(4,.5)(0,4.5)
\psline(4,.5)(5,1.5)
\psline(3,1.5)(4,2.5)
\psline(2,2.5)(3,3.5)
\psline(1,3.5)(2,4.5)
\end{pspicture}}

\newcommand\treeBACD{
\setlength{\unitlength}{3pt}
\psset{unit=3pt}
\psset{runit=2pt}
\psset{linewidth=0.2}
\begin{pspicture}(0,-.5)(6,5)
\psline(4,-1)(4,.5)
\psline(4,.5)(1,3.5)
\psline(4,.5)(5,1.5)
\psline(3,1.5)(4,2.5)
\psline(2,2.5)(4,4.5)
\psline(3,3.5)(2,4.5)
\end{pspicture}}

\newcommand\treeABDA{
\setlength{\unitlength}{3pt}
\psset{unit=3pt}
\psset{runit=2pt}
\psset{linewidth=0.2}
\begin{pspicture}(0,-.5)(7,4.5)
\psline(4,-1)(4,.5)
\psline(4,.5)(.5,4)
\psline(4,.5)(6.5,3)
\psline(2.5,2)(3.5,3)
\psline(1.5,3)(2.5,4)
\psline(5.5,2)(4.5,3)
\end{pspicture}}

\newcommand\treeBADA{
\setlength{\unitlength}{3pt}
\psset{unit=3pt}
\psset{runit=2pt}
\psset{linewidth=0.2}
\begin{pspicture}(0,-.5)(6,4.5)
\psline(3,-1)(3,.5)
\psline(3,.5)(.5,3)
\psline(3,.5)(5.5,3)
\psline(1.5,2)(3.5,4)
\psline(2.5,3)(1.5,4)
\psline(4.5,2)(3.5,3)
\end{pspicture}}

\newcommand\treeADAB{
\setlength{\unitlength}{3pt}
\psset{unit=3pt}
\psset{runit=2pt}
\psset{linewidth=0.2}
\begin{pspicture}(0,-.5)(6,4.5)
\psline(3,-1)(3,.5)
\psline(3,.5)(.5,3)
\psline(3,.5)(5.5,3)
\psline(1.5,2)(2.5,3)
\psline(4.5,2)(2.5,4)
\psline(3.5,3)(4.5,4)
\end{pspicture}}

\newcommand\treeADBA{
\setlength{\unitlength}{3pt}
\psset{unit=3pt}
\psset{runit=2pt}
\psset{linewidth=0.2}
\begin{pspicture}(0,-.5)(6.5,4.5)
\psline(3,-1)(3,.5)
\psline(3,.5)(.5,3)
\psline(3,.5)(6.5,4)
\psline(1.5,2)(2.5,3)
\psline(4.5,2)(3.5,3)
\psline(5.5,3)(4.5,4)
\end{pspicture}}

\newcommand\treeDCAB{
\setlength{\unitlength}{3pt}
\psset{unit=3pt}
\psset{runit=2pt}
\psset{linewidth=0.2}
\begin{pspicture}(0,-.5)(5.5,5)
\psline(2,-1)(2,.5)
\psline(2,.5)(1,1.5)
\psline(2,.5)(5,3.5)
\psline(3,1.5)(2,2.5)
\psline(4,2.5)(2,4.5)
\psline(3,3.5)(4,4.5)
\end{pspicture}}

\newcommand\treeDCBA{
\setlength{\unitlength}{3pt}
\psset{unit=3pt}
\psset{runit=2pt}
\psset{linewidth=0.2}
\begin{pspicture}(0,-.5)(6,5)
\psline(2,-1)(2,.5)
\psline(2,.5)(1,1.5)
\psline(2,.5)(6,4.5)
\psline(3,1.5)(2,2.5)
\psline(4,2.5)(3,3.5)
\psline(5,3.5)(4,4.5)
\end{pspicture}}

\newcommand\combL{
\setlength{\unitlength}{3pt}
\psset{unit=3pt}
\psset{runit=2pt}
\psset{linewidth=.2}
\begin{pspicture}(-1,-.5)(7,8)
\psline(5,-.5)(5,1)
\psline(5,1)(0,6)
\psline(5,1)(6,2)
\put(6.5,4){.}
\put(5.5,5){.}
\put(4.5,6){.}
\psline(2,4)(3,5)
\psline(1,5)(2,6)
\end{pspicture}}

\newcommand\combR{
\setlength{\unitlength}{3pt}
\psset{unit=3pt}
\psset{runit=2pt}
\psset{linewidth=.2}
\begin{pspicture}(-1,-.5)(7,8)
\psline(1,-.5)(1,1)
\psline(1,1)(6,6)
\psline(1,1)(0,2)
\psline(2,2)(1,3)
\put(2.5,5.5){.}
\put(3.5,6.5){.}
\put(4.5,7.5){.}
\psline(5,5)(4,6)
\end{pspicture}}

\newcommand{\lgraft}[2]
{\setlength{\unitlength}{4pt}
\psset{unit=5pt}
\psset{runit=4pt}
\psset{linewidth=.1}
\begin{pspicture}(0,0)(4,4)
\psline(2.1,.9)(1.4,1.6)
\put(2.6,0){$#1$}
\put(0.5,2.6){$#2$}
\end{pspicture}}

\newcommand{\rgraft}[2]
{\setlength{\unitlength}{4pt}
\psset{unit=5pt}
\psset{runit=4pt}
\psset{linewidth=.1}
\begin{pspicture}(0,0)(4,4)
\psline(1.5,.9)(2.2,1.6)
\put(.4,0){$#1$}
\put(3,2.6){$#2$}
\end{pspicture}}

\newcommand{\lvertexgraft}[1]
{\setlength{\unitlength}{3pt}
\psset{unit=3pt}
\psset{runit=2pt}
\psset{linewidth=0.2}
\begin{pspicture}(0,-.5)(4.5,4.5)
\psline(3,-1)(3,.5)
\psline(3,.5)(2,1.5)
\psline(3,.5)(4,1.5)
\put(1,3.2){\footnotesize$#1$}
\end{pspicture}}

\newcommand{\rvertexgraft}[1]
{\setlength{\unitlength}{3pt}
\psset{unit=3pt}
\psset{runit=2pt}
\psset{linewidth=0.2}
\begin{pspicture}(0,-.5)(4.5,4.5)
\psline(1,-1)(1,.5)
\psline(1,.5)(0,1.5)
\psline(1,.5)(2,1.5)
\put(2.5,3.2){\footnotesize$#1$}
\end{pspicture}}

\newcommand{\rvertexgraftsmall}[1]
{\setlength{\unitlength}{2pt}
\psset{unit=2pt}
\psset{runit=1.2pt}
\psset{linewidth=0.2}
\begin{pspicture}(0,-.5)(5,3.5)
\psline(1,-1)(1,.5)
\psline(1,.5)(0,1.5)
\psline(1,.5)(2,1.5)
\put(2.7,3.2){\tiny$#1$}
\end{pspicture}}

\newcommand{\combLgraft}[3]
{\setlength{\unitlength}{5pt}
\psset{unit=6pt}
\psset{runit=5pt}
\psset{linewidth=.1}
\begin{pspicture}(-1,-.5)(7,8)
\psline(5,-.5)(5,1)
\psline(5,1)(0,6)
\psline(5,1)(6,2)
\put(5.8,3.3){.}
\put(5,4.1){.}
\put(4.1,4.9){.}
\psline(2,4)(3,5)
\psline(1,5)(2,6)
\put(2.4,7.8){\footnotesize$#1$}
\put(3.8,6.3){\footnotesize$#2$}
\put(7.3,3){\footnotesize$#3$}
\end{pspicture}}

\newcommand{\combRgraft}[3]
{\setlength{\unitlength}{5pt}
\psset{unit=6pt}
\psset{runit=5pt}
\psset{linewidth=.1}
\begin{pspicture}(-1,-.5)(9,7)
\psline(3,-.5)(3,1)
\psline(3,1)(8,6)
\psline(3,1)(2,2)
\psline(4,2)(3,3)
\psline(7,5)(6,6)
\put(4.4,4.4){.}
\put(5.2,5.2){.}
\put(6,6){.}
\put(0.2,2.9){\footnotesize$#1$}
\put(1.7,4.6){\footnotesize$#2$}
\put(5.4,7.9){\footnotesize$#3$}
\end{pspicture}}

\begin{document}

\title{Groups of tree-expanded series}
 
\author{
Alessandra Frabetti \\
Universit\'e de Lyon,
Universit\'e Lyon1, CNRS, \\ 
UMR 5208 Institut Camille Jordan, \\
B\^atiment du Doyen Jean Braconnier, \\ 
43, blvd du 11 novembre 1918,
F-69622 Villeurbanne Cedex,
France \\ 
{\small \tt frabetti@math.univ-lyon1.fr}}
\date{\today}
 
\maketitle


\begin{abstract}
  In \cite{BFqedren,BFqedtree} we introduced three Hopf algebras on
  planar binary trees related to the renormalization of quantum
  electrodynamics. One of them, the algebra $\Ha$, is commutative,
  and is therefore the ring of coordinate functions of a proalgebraic 
  group $\Ga$. The other two algebras, $\He$ and $\Hp$, are free 
  non-commutative. Therefore their abelian quotients are the coordinate 
  rings of two proalgebraic groups $\Ge$ and $\Gp$. In this paper we describe 
  explicitly these groups.

  Using two monoidal structures and a set-operad structure on planar
  binary trees, we show that these groups can be realized on formal
  series expanded over trees, and that the group laws are
  generalizations of the multiplication and the composition of usual
  series in one variable. Therefore we obtain some new groups of invertible
  tree-expanded series and of tree-expanded formal diffeomorphisms
  respectively.

  The Hopf algebra describing the renormalization of the electric charge 
  corresponds to the subgroup of tree-expanded formal diffeomorphisms
  formed of the translations, which fix the zero, by some particular 
  tree-expanded series which remind the proper correlation functions 
  in quantum field theory. 
  In turn, the group of tree-expanded formal diffeomorphisms and some 
  of its subgroups give rise to new Hopf algebras on trees.

  All the constructions are done in a general operad-theoretic setting, 
  and then applied to the specific duplicial operad on trees.
\end{abstract}

\bigskip 

\bigskip

\bigskip 


\section*{Introduction}

In \cite{BFqedtree}, C.~Brouder and the author introduced three Hopf 
algebras $\He$, $\Hp$ and $\Ha$ related to the renormalization of 
perturbative quantum electrodynamics. They are constructed on planar 
binary trees, where each tree represents a suitable sum of QED Feynman 
graphs. 
The abelian quotients $\He_{ab}$, $\Hp_{ab}$, and the algebra $\Ha$, 
are commutative Hopf algebras, and therefore they are coordinate rings 
for some proalgebraic groups, that we denote respectively by 
$\Ge$, $\Gp$ and $\Ga$. 
By proalgebraic group, we mean a functor $G$ which associates a group $G(A)$ 
to any associative, unital and commutative algebra $A$, and which is 
represented by a commutative Hopf algebra ${\cal{H}}$ which is not necessarily 
finitely generated. In other words, the group $G(A)$ is isomorphic 
to the group $\Hom_{Alg}({\cal{H}},A)$ of algebra homomorphisms, 
considered with the convolution product inherited from the coalgebra 
structure of ${\cal{H}}$, cf. \cite{Abe,Borel}. 

In this paper, we describe explicitly these groups as group functors, 
and show that they can be considered as generalizations of the groups 
of formal series in one variable, endowed with the multiplication and 
the composition. 
To do this, we introduce the notion of tree-expanded series. 

Let $A$ be an associative, unital and commutative algebra over the field 
$\Q$. 
Denote by $Y$ the set of all trees, and by $A[[Y]]$ the vector space of 
sequences $(a_t)_{t\in Y}$ where $a_t \in A$. For our purpose, it is 
convenient to write a sequence $(a_t)$ as a formal series 
$a(x)=\sum_{t\in Y} a_t\ x^t$, where $x$ is a formal variable. 
Here the word ``series'' is an abuse of language, since $x^t$ is just a symbol,
for instance $x^{\!\!\!\!\treeAB}$.
We call such sequences {\em tree-expanded series\/} in a variable $x$, 
with coefficients in $A$. 
These series behave like usual ones. We can multiply them, and compose those 
which have zero constant term. Here, the constant term is the coefficient 
of the symbol $x^{\treeO}$, where $\treeO$ is the tree given by the sole root. 

As affine varieties, we can identify the groups as follows. 
If we denote by $\Y$ the set of trees different from $\treeO$, the groups 
$\Ge(A)$ and $\Gp(A)$ are both isomorphic to the subset $x^{\treeO}+A[[\Y]]$. 
On the other side, if we denote by $\YY$ the set of trees different 
from $\treeO$ and $\treeA$, the group $\Ga(A)$ is isomorphic 
to the subset of $x^{\treeA}+A[[\YY]]$ containing the tree-expanded 
series of the form 
\begin{align*}
\varphi(x)=x^{\treeA}+\sum_{t\in \Y} \varphi_t\ x^{\!\!\lvertexgraft{t}}, 
\end{align*}
where the coefficients satisfy the condition
\begin{align*}
\varphi_t = \varphi_{\rvertexgraftsmall{t_1}}\ 
\varphi_{\rvertexgraftsmall{t_2}}\ 
\cdots\varphi_{\rvertexgraftsmall{t_n}}
\qquad\mbox{if\quad $t=\combLgraft{t_1}{t_2}{t_n}$}. 
\end{align*}

In this paper we describe the group laws of $\Ge(A)$, $\Gp(A)$ and $\Ga(A)$ 
presented as sets of tree-expanded series, and their relationship with
their analogue groups of usual formal series. 

Seen as representable group functors, the isomorphisms
$G(A)\cong\Hom_{Alg}({\cal{H}},A)$ bring a tree-expanded formal 
series $f(x)=\sum_{t\in Y} f_t\ x^t$ into the algebra homomorphism from 
${\cal{H}}$ to $A$ which associates to the tree $t\in {\cal{H}}$ 
the coefficient $f_t\in A$. 
\bigskip

This situation is in fact not peculiar to trees. 
We show that the same kind of constructions can be done on the set 
$A[[\P]]$ of $\P$-expanded series, if $\P$ is a graded collection of 
finite sets with suitable properties. 
These turn out to be: a set-operad structure on $\P$, and two associative 
binary operations in $\P(2)$. The results on trees are then obtained 
by considering the duplicial operad generated by the {\em over} and 
{\em under} grafting operations on trees. 
Similar results could be obtained for other operads, and in particular 
for the diassociative operad describing dialgebras, however we do not 
investigate the resulting groups of series in this paper. 
\bigskip 

One of the key resulting groups, that of formal diffeomorphisms, 
requires only an operad structure. For algebraic operads, essentially 
the same construction was considered by F.~Chapoton in \cite{Chapoton} 
and by P.~van~der~Laan in \cite{VanDerLaan} with different motivations 
and applications. 
Chapoton specifies to the {\em pre-Lie operad\/} of rooted trees, and uses 
the rooted tree-expanded series to describe the flow of a linear vector field 
on an affine space. Van~der~Laan describes the non-symmetric case and  
introduces the non-commutative lift of the coordinate rings. 
We will comment briefly the relationship between the different constructions 
in Section~\ref{set-operad}. 
\bigskip\bigskip  

\noindent
{\bf Notations.} 
All vector spaces and algebras are defined over the field $\Q$ of 
rational numbers, although the algebras $\He$, $\Hp$ and $\Ha$ 
were originally defined over the field of complex numbers. In fact, 
this restriction was not necessary. 

For any set $X$, we denote by $\Q X$ the vector space spanned by $X$, by 
$\Q\langle X \rangle$ the tensor algebra on $X$ (non-commutative polynomials), 
and by $\Q[X]$ the symmetric algebra on $X$ (commutative polynomials). 
\bigskip

\noindent
{\bf Acknowledgments.} 
The author is deeply indebted to Olivier Mathieu and Fr\'ederic Chapoton 
for their careful reading and useful suggestions on the paper. 
She is particularly grateful to Fr\'ederic Chapoton for his interest on 
the subject, and his computations of some peculiar tree-expanded series 
which will be presented in a separate contribution, 
cf.~\cite{ChapotonFrabetti}. 

The author warmly thanks the Swiss National Foundation for Scientific 
Research for the support of her visit to the Mathematics Department 
of Lausanne University in 2001, where this work originated. She thanks 
as well the members of the Mathematics Department of Lausanne University 
for their hospitality. 


\tableofcontents


\section{Group of invertible tree-expanded series}

In this section, we recall how to associate a group of series to
a graded monoid, and discuss the relationship with the group of invertible 
usual series in one variable. 

Then we consider two graded monoids built on trees, with the {\em over} 
and {\em under} grafting operations, and describe the resulting groups 
of tree-expanded series. We show that the coordinate rings of these 
groups coincide with the ``QED propagator'' Hopf algebras introduced 
in \cite{BFqedtree}.

\subsection{Group of series expanded over a graded monoid}
\label{monoid}

Let $\M(n)$ be a collection of finite sets indexed by non-negative integers, 
and set $\M=\sqcup_{n\geq 0} \M(n)$. 
Then $\M$ is a graded monoid if it is endowed with an associative graded 
operation $\cdot : \M\times\M\longrightarrow\M$ and a neutral element 
$e\in\M(0)$. 

Let $A$ be an associative, commutative and unital algebra over $\Q$, and 
denote by $A[[\M]]$ the vector space of $\M$-expanded series
$f(x)=\sum_{p\in\M} f_p\ x^p$, with $f_p\in A$.
We define a multiplication in $A[[\M]]$ by setting 
\begin{align}
\label{def-multiplication}
(f\cdot g)(x) &:= \sum_{p\in\M}\ \sum_{q\in\M} f_p\ g_q\ x^{p\cdot q}.  
\end{align}
This series is well defined because for a given $u\in\M$ there exist finitely 
many elements $p,q\in\M$ such that $p\cdot q=u$. 
Then $A[[\M]]$ forms a unital algebra, with unit $x^e$. 

\begin{theorem} {\bf (obvious)} 
\label{Gi_M}
Set $\MM= \sqcup_{n\geq 1} \M(n)$. 
The subset $\Gi_{\M}(A):=x^e+A[[\MM]]$ forms a group.
\end{theorem}

We call $\Gi_{\M}(A)$ the group of {\em invertible $\M$-series\/}, with 
coefficients in $A$. 
\bigskip 

\noindent{\bf Example.} 
The first example of such a group is constructed from the graded monoid 
$\N$ of non-negative integers, with $\N(n)=\{n\}$ and the 
addition. We then obtain the group $\Gi(A)=1+xA[[x]]$ of usual 
invertible formal series. 
\bigskip

The construction of the group $\Gi_{\M}(A)$ is functorial in $\M$. 
Since the grading $\pi:\M\longrightarrow\N$ is a canonical 
morphism of graded monoids, and any $p\in\M(1)$ gives a section 
$i_p:\N\longrightarrow\M$ by $n\mapsto p^n$, 
we can compare the series expanded over $\M$ with the usual formal series. 

\begin{proposition}
\label{Gi_M->Gi}
For any graded monoid $\M$, there is a canonical morphism of groups 
$\pi:\Gi_{\M}(A)\longrightarrow\Gi(A)$. 
Moreover, any element $p\in\M(1)$ gives a section  
$i_p:\Gi(A)\longrightarrow\Gi_{\M}(A)$. 
\end{proposition}
\bigskip 

The construction of the group $\Gi_{\M}(A)$ is also functorial in $A$. 
Indeed, the group functor $\Gi_{\M}$ can be represented by a commutative 
Hopf algebra $\Hi_{\M}$, i.e., for any associative, commutative and 
unital algebra $A$, we have 
\begin{align*}
\Gi_{\M}(A) &\cong \Hom_{Alg}(\Hi_{\M},A).  
\end{align*}
As an algebra, $\Hi_{\M}$ is a polynomial algebra with generators indexed by 
$\MM$. It is remarkable that it admits a straightforward lift to a 
non-commutative polynomial algebra. 

Given a commutative Hopf algebra $\H$, by {\em non-commutative lift\/} 
of $\H$ we mean a non-commutative Hopf algebra $\Hnc$ such that $\H$ is 
the abelian quotient of $\Hnc$. 
The {\em abelian quotient\/} of a Hopf algebra $\Hnc$ is the commutative 
algebra $\Hnc_{ab}$ obtained as the quotient of $\Hnc$ on the ideal 
generated by the commutator $[\Hnc,\Hnc]$. It is automatically a Hopf algebra. 

In the case of $\Hi_{\M}$, we consider the free associative algebra 
$\Hinc_{\M}=\Q\langle\MM\rangle$, with generators graded by the grading 
of $\M$. The neutral element $e\in\M(0)$ is identified with the unit $1$. 
Define a coproduct 
$\Di_{\M}:\Hinc_{\M}\longrightarrow\Hinc_{\M}\otimes\Hinc_{\M}$ 
as the algebra morphism given on the generators $u$ by 
\begin{align*}
\Di_{\M}(u) &:= \sum_{u=p\cdot q} p \otimes q. 
\end{align*}
Define also a counit $\varepsilon:\Hinc_{\M}\longrightarrow\Q$ as the 
algebra morphism given on all $u\neq e$ by $\varepsilon(u)=0$. 

\begin{proposition}
\label{Hi_M}
The algebra $\Hinc_{\M}$ is a (non-commutative) graded and connected Hopf 
algebra. It is cocommutative if and only if the monoid is abelian. 

The coordinate ring of $\Gi_{\M}$ is the abelian quotient $\Hi_{\M}=\Q[\MM]$ 
of the Hopf algebra $\Hinc_{\M}$. 
\end{proposition}

\noindent{\bf Example.} 
In particular, the coordinate ring of the group functor $\Gi$ is the
polynomial algebra $\Hi=\Q[b_1,b_2,\dots]$ on one generator $b_n$ in
each degree $n\geq 1$. The coproduct on the generators of $\Hi$ is
\begin{align*}
\Di b_n = \sum_{k=0}^n b_k \otimes b_{n-k}, \qquad (b_0=1),
\end{align*}
and the counit is $\varepsilon(b_n)=0$ for $n\geq 1$.  This Hopf
algebra is the unique free commutative and cocommutative
Hopf algebra on infinitely many generators spanned by natural numbers,
commonly known as the Hopf algebra of symmetric functions. It
is well known to admit a non-commutative lift to the tensor algebra
$\Hinc=\Q\langle b_1,b_2,\dots\rangle$ which is still cocommutative.
\bigskip

Since the correspondance between proalgebraic groups and their representative 
Hopf algebras is contravariant, the relationship between the Hopf algebras 
$\Hinc_{\M}$ and $\Hinc$ can be find by reversing the morphisms of 
Proposition~\ref{Gi_M->Gi}. 

\begin{proposition}
\label{Hi->Hi_M}
There is a canonical morphism of Hopf algebras 
$\Hinc\longrightarrow\Hinc_{\M}$. 
Moreover, each element $p\in\M(1)$ determines a section 
$\Hinc_{\M}\longrightarrow\Hinc$. 
\end{proposition}

\begin{proof}
The morphism $\Hinc\longrightarrow\Hinc_{\M}$ is explicitly given by 
\begin{align*}
b_n \longmapsto \sum_{p\in\M(n)} p, 
\end{align*} 
and it is injective if there is an element $p\in\M(1)$. 
Then, its section $\Hinc_{\M}\longrightarrow\Hinc$ 
is given by the map 
\begin{align*}
u\in \M(n) \longmapsto 
\left\{\begin{array}{ll} b_n & \mbox{if $u=p^n$} \\ 
0 & \mbox{otherwise}
\end{array}\right. .
\end{align*}
The verification that these two maps are morphisms of Hopf algebras is 
straightforward. 
\end{proof}


\subsection{Graded monoid of trees and invertible tree-expanded 
series}
\label{tree-invertible}

In this paper, by a ``tree'' we will always mean a {\em planar binary 
rooted tree\/}, that is, a connected planar graph without loops, having 
internal vertices of valence 3 and a preferred external vertex called 
the {\em root\/}. For any tree $t$, we call {\em order of $t$\/} the number 
$|t|$ of its internal vertices. For any $n\geq 0$, we then denote by $Y_n$ 
the set of trees of order $n$. In particular, there is only one tree 
$\treeO$ with no internal vertex, called the {\em root tree\/}, and there 
is only one tree $\treeA$ with 1 internal vertex, called the {\em vertex 
tree\/}. The sets $Y_2$ and $Y_3$ are the following:
\begin{align*}
Y_2 &= \left\{ \treeAB, \treeBA \right\}, \\ 
Y_3 &= \left\{ \treeABC, \treeBAC, \treeACA, \treeCAB, \treeCBA \right\}. 
\end{align*}

On the set $Y=\cup_{n\geq 0} Y_n$ of all trees, let us consider the over 
and under products introduced by J.-L.~Loday in~\cite[Sec.~1.5]{LodayATree}. 
We recall that for any planar binary trees $t$ and $s$, the tree 
$t$~{\em over}~$s$ is the grafting
\begin{align*}
t \over s = \lgraft{s}{t}
\end{align*}
of the root of $t$ on the left-most leaf of $s$, while the tree
$t$~{\em under}~$s$ is the grafting
\begin{align*}
t \under s = \rgraft{t}{s}
\end{align*}
of the root of $s$ on the right-most leaf of $t$. 
Of course, the operations over and under are isomorphic, 
and the isomorphism is given by the axial symmetry of the trees along 
their roots. 

The over and under products on trees are associative, non-commutative, 
with unit given by the root tree $\treeO$. Moreover they are graded with 
respect to the order of trees, in the sense that
\begin{align}
\label{order}
|t \over s| = |t|+|s| \quad\mbox{and}\quad |t\under s|=|t|+|s|. 
\end{align} 
Therefore $(Y,\over)$ and $(Y,\under)$ are two graded monoids.
Applying Theorem~\ref{Gi_M}, we then obtain two non-abelian group laws
on the set of tree-expanded series
\begin{align*}
\Gi_Y(A) &:= \{ f(x)=\sum_{t\in Y} f_t\ x^t,\ f_t\in A,\ f_{\treeO}=1 \}, 
\end{align*}
with multiplications denoted by $\over$ and $\under$, and with unit 
$x^{\treeO}$. Denote these groups as 
\begin{align*}
\Go(A):=(\Gi_Y(A),\over) \quad\mbox{and}\quad
\Gu(A):=(\Gi_Y(A),\under). 
\end{align*}

\begin{theorem}
The non-commutative lift of the coordinate rings of the two group  
functors $\Go$ and $\Gu$ coincide with the QED propagator Hopf algebras 
$\Hp$ and $\He$.
\end{theorem}

In other words, if we denote by $\Gp$ and $\Ge$ the group functors 
represented respectively by the Hopf algebras $\Hp_{ab}$ and $\He_{ab}$, 
then the groups are exactly $\Gp=\Go$ and $\Ge=\Gu$. 
\bigskip 

\begin{proof}
Let us recall, from \cite{BFqedtree}, the definition of the ``QED propagator 
Hopf algebras'' $\Hp$ and $\He$ on planar binary trees. As algebras, they
are both isomorphic to the free non-commutative algebra generated by
all trees different from $\treeO$, that is $\Hp=\He=\Q\langle\Y\rangle$, 
where $\Y=\cup_{n\geq 1} Y_n$. 
Since $Y_0=\{\treeO\}$, we also identify the root-tree $\treeO$ to the
unit $1$ and write $\Hp=\He \cong \Q\langle Y \rangle/(\treeO-1)$.

The coalgebra structures are given by the ``pruning coproducts''
$\Dip:\Hp\longrightarrow\Hp\otimes\Hp$ and
$\Die:\He\longrightarrow\He\otimes\He$ defined as the dual operations
respectively to the over and under products of trees. That is, 
they are defined on the generators $u\in \Y$ by
\begin{align}
\label{Di^pe} 
\Dip(u) = \sum_{t\over s=u} t \otimes s \qquad\mbox{and}\qquad 
\Die(u) = \sum_{t\under s=u} t \otimes s. 
\end{align}
The counits $\varepsilon:\Hp\longrightarrow \Q$ and
$\varepsilon:\He\longrightarrow \Q$ are dual to the unit $\treeO$,
that is $\varepsilon(\treeO)=1$ and $\varepsilon(t)=0$ if $t\neq
\treeO$, and the antipodes are then defined in a standard recursive
way.

Comparing the pruning coproducts in $\Hp$ and $\He$ given by 
Eq.~(\ref{Di^pe}) with the over and under multiplications in $\Gp(A)$ 
and $\Ge(A)$ obtained from Eq.~(\ref{def-multiplication}) for the 
{\em over\/} and {\em under\/} products, it is clear that the abelian 
quotients $\Hp_{ab}$ and $\He_{ab}$ of the two non-commutative Hopf algebras 
$\Hp$ and $\He$ are respectively the coordinate rings of the two group 
functors $\Go$ and $\Gu$.
\end{proof}

Let us call {\em order map\/} the map $|\ |: Y\longrightarrow \N$
which sends each tree $t\in Y$ to its order $|t|\in\N$.
Formulas~(\ref{order}) say that it is a morphism of graded monoids,
and in fact it coincides with the projection $\pi$ of
Proposition~\ref{Gi_M->Gi}.  Since $Y_1=\{\treeA\}$, we conclude that
the order map induces two surjective morphisms of groups from $\Go(A)$ 
and $\Gu(A)$ respectively to $\Gi(A)$.

The vertex tree then determines a section for each of the two
projections, namely the maps $\io,\iu:\Gi(A) \longrightarrow \Gi_Y(A)$ 
such that 
\begin{align*}
\io(x^n) = x^{\treeA^{\over n}} = x^{\combL} \qquad\mbox{and}\qquad  
\iu(x^n) = x^{\treeA^{\under n}} =x^{\combR}. 
\end{align*} 
Let us call these trees respectively the left and the right {\em comb 
trees\/}. 
Therefore the inclusions $\iu$ and $\io$ identify the group $\Gi(A)$ of
invertible series with the two subgroups of $\Gu(A)$ and $\Go(A)$ made of
tree-expanded series expanded only on the comb trees.

The result of Proposition~\ref{Hi->Hi_M} on the dual Hopf algebras $\He$
and $\Hp$, namely that they contain $\Hinc$ as a Hopf subalgebra, was
already shown in \cite{BFK}.


\section{Tree-expanded formal diffeomorphisms}

In this section we introduce a group of series associated to any 
non-symmetric operad in the category of sets, and discuss the relationship 
with the group of usual formal diffeomorphisms in one variable. 

Then we apply the results to the duplicial operad, generated by the 
{\em over} and {\em under} operations on trees. The resulting group of 
tree-expanded diffeomorphisms is bigger than the group $\Ga$ that 
we are looking for. However we describe it explicitly, because its 
coordinate ring admits a non-commutative lift which gives rise to 
a new graded and connected Hopf algebra on trees, which is neither 
commutative nor cocommutative. 

\subsection{Group of series expanded over a set-operad}
\label{set-operad}

Let us recall the definition of a (non-symmetric) operad in the category of 
sets. 
Let $\P(n)$ be a collection of finite sets, indexed by positive integers, 
and set $\P=\sqcup_{n\geq 1} \P(n)$. 
Then $\P$ is a set-operad if there exist compositions maps 
\begin{align*}
\gamma \equiv \gamma_{n;m_1,\dots,m_n}:\ 
& \P(n)\times\P(m_1)\times\dots\times\P(m_n) 
\longrightarrow \P(m_1+\dots +m_n) \\
&(p;q_1,\dots,q_n) \mapsto \gamma(p;q_1,\dots,q_n)
\end{align*}
satisfying the associative condition 
\begin{multline*}
 \gamma(\gamma(p;q_1,\dots,q_n) ; u^1_1,\dots,u^1_{m_1},u^2_1,\dots,u^2_{m_2},
\dots,u^n_1,\dots,u^n_{m_n})
\\= \gamma(p ; \gamma(q_1;u^1_1,\dots,u^1_{m_1}), \dots ,
\gamma(q_n;u^n_1,\dots,u^n_{m_n}))
\end{multline*}
and an operation $\id\in\P(1)$ acting as the identity, that is 
\begin{align*}
  \gamma(\id;q)=q \quad \text{and} \quad
  \gamma(p;\id,\dots,\id)=p.
\end{align*}
For any $p\in\P(n)$, set $|p|=n$ and call it the {\em order of $p$}. 

The canonical example of a set-operad is the Endomorphism operad $\Endom_S$ 
of a set $S$, given by the collection of set-maps $\Endom_S(n)=\Hom(S^n,S)$, 
together with the usual composition of maps. 
If $\P$ is a set-operad, a morphism from $\P$ to the endomorphism operad 
$\Endom_S$ defines on the set $S$ the structure of a $\P$-monoid. 
\bigskip 

Let $A$ be an associative, commutative and unital algebra over $\Q$, 
and denote by $A[[\P]]$ the vector space of $\P$-expanded series 
$\varphi(x)=\sum_{p\in\P}f_p\ x^p$, with $f_p\in A$. 
We define the {\em composition\/} in $A[[\P]]$ by 
\begin{align}
\label{P-composition}
(\varphi\circ\psi)(x) 
&:= \sum_{p\in\P}\ \ \sum_{q_1,\dots,q_{|p|}\in\P}
\varphi_{p} \psi_{q_1}\dots \psi_{q_n}\ x^{\gamma(p;q_1,\dots,q_{|p|})}. 
\end{align}
 
\begin{theorem}
\label{Gd_P}
Set $\PP=\sqcup_{n\geq 2} \P(n)$. 
The subset $\Gd_{\P}(A):=x^{\id}+A[[\PP]]$ forms a group, with unit $x^{\id}$. 
Moreover the composition respects the shifted grading $\deg(p)=|p|-1$.
\end{theorem}

We call $\Gd_{\P}(A)$ the group of {\em $\P$-expanded formal 
diffeomorphisms\/}.
\bigskip 

\begin{proof}
The series $x^{\id}\in\Gd_{\P}(A)$ is obviously a unit. 
The subset $\Gd_{\P}(A)$ is obviously closed under $\circ$. 
To show that it is a group, it suffices to show that the composition in 
$A[[\P]]$ is associative and respects the shifted grading. 
In fact, its series $\varphi$ have invertible constant term 
$\varphi_{\id}=1$, and therefore their compositional inverse can 
be found recursively. 

The associativity of the composition $\circ$ is ensured by the associativity 
of the operadic composition $\gamma$. 
In fact, given three $\P$-expanded series $\varphi,\psi,\eta$, we have 
\begin{align*}
[(\varphi\circ\psi)\circ\eta](x) &= 
\underset{u_1,\dots,u_m \in\P}{\underset{q_1,\dots,q_n \in\P}{\sum_{p\in\P}}} 
\varphi_p\ \psi_{q_1}\dots\psi_{q_n}\ \eta_{u_1} \dots \eta_{u_m}\ 
x^{\gamma(\gamma(p;q_1,\dots,q_n);u_1,u_2,\dots,u_m)}, 
\end{align*}
where $n=|p|$ and $m=|\gamma(p;q_1,\dots,q_n)|=\sum_{i=1}^n |q_i|$. 
On the other side, we have 
\begin{align*}
[\varphi\circ(\psi\circ\eta)](x) &= 
\underset{q_1,\dots,q_n\in\P}{\sum_{p\in\P}}
\underset{v^n_1,\dots,v^n_{|q_n|}\in\P}{\underset{\dots..\quad}
{\sum_{v^1_1,\dots,v^1_{|q_1|}\in\P}}} 
\varphi_p\ \psi_{q_1}\dots\psi_{q_n}\ \eta_{v^1_1}\dots\eta_{v^1_{|q_1|}} 
\dots\eta_{v^n_1}\dots\eta_{v^n_{|q_n|}} \times \\ 
& \hspace{5cm} \times\ x^{\gamma(p;\gamma(q_1;v^1_1,\dots,v^1_{|q_1|}),\dots,
\gamma(q_n;v^n_1,\dots,v^n_{|q_n|}))}, 
\end{align*}
where the total number of operations $v^i_j$ is $|q_1|+\dots+|q_n|=m$. 
If we call $v^1_{k_1}:= u_{k_1}$, and 
$v^i_{k_i}:= u_{|q_1|+\dots+|q_{i-1}|+k_i}$, for $i=2,\dots,n$, 
all the factors $\eta_{v^i_{k_i}}$ of the second term have a corresponding 
factor $\eta_{u_j}$ in the first term, and the associativity of the 
composition $\gamma$ guarantees that the exponents coincide. 

The shifted grading ensures that the composition $\circ$ is graded. 
In fact, for any operations $p$ and $q_1,\dots,q_{|p|}$, we have 
\begin{align*}
\deg(p) + \deg(q_1 \dots q_{|p|}) 
&= (|p|-1) + \sum_{i=1}^{|p|} (|q_i|-1) 
= \sum_{i=1}^{|p|} |q_i|-1 \\  
&= \deg(\gamma(p;q_1,\dots,q_{|p|})).   
\end{align*}
\end{proof}

\noindent{\bf Example.} 
The simplest example of such a group is that of usual formal diffeomorphisms. 
To see how it arises from a set-operad, consider the set-operad 
$\N^*=\bigsqcup_{n\geq 1} \{n\}$ of positive integers, with the 
compositions $\gamma(n;m_1,\dots,m_n)=m_1+\dots+m_n$ and identity $\id=1$.  
It is the quadratic set-operad generated by one associative operation 
in degree $2$. A $\N^*$-monoid is a usual associative monoid. 
From a set-operad $\P$ one can define an algebraic operad $\Q\P$ by 
considering the collection of $\Q$-vector spaces $\Q\P(n)$ spanned by 
the finite sets, and extending linearly the compositions $\gamma$ to maps
$\Q\P(n)\otimes\Q\P(m_1)\otimes\dots\otimes\Q\P(m_n) \longrightarrow
\Q\P(m_1+\dots +m_n)$.
The algebraic operad associated to the set-operad $\N^*$ is the
operad $\As$ of associative algebras. Therefore we will use the symbol 
$\As$ also to denote the set-operad. 

For $\P=\As$ we have $A[[\As]]=xA[[x]]$ and $A[[\bar{\As}]]=x^2A[[x]]$. 
Therefore the group $\Gd_{\As}(A)$ is the group $\Gd(A)= x + x^2 A[[x]]$ 
of usual formal diffeomorphisms 
$\varphi(x)=x+\sum_{n\geq 2}\varphi_n\ x^n$ (tangent to the identity), 
with coefficients in $A$, considered with the composition (or substitution) 
\begin{align}
\label{composition}
(\varphi\circ\psi)(x) &= \varphi(\psi(x)) 
= \sum_{n=1}^\infty \varphi_n \psi(x)^n 
= \sum_{n=1}^\infty 
\left( \sum_{m=1}^n \underset{k_1,\dots,k_m\geq 1}{\sum_{k_1+\dots+k_m=n}} 
\varphi_m\ \psi_{k_1}\dots \psi_{k_m} \right)\ x^n, 
\end{align}
and with unit given by the series $x$. 
\medskip 

Another example can be constructed from the diassociative set-operad $\Dias$, 
whose algebraic extension was introduced by J.-L.~Loday in 
\cite{LodayDialgebras} and gives rise to dialgebras. 
It is given by the collection of sets $\Dias(n)=\{1,\dots,n\}$, 
with compositions 
\begin{align*}
\Dias(n)\times\Dias(m_1)\times\cdots\times\Dias(m_n)
& \longrightarrow \Dias(m_1+\cdots +m_n) \\ 
(i;j_1,...,j_n) & \longmapsto j_i. 
\end{align*}
However we do not investigate here the associated group of formal 
diffeomorphisms. A non-trivial example of such a group is explained in 
details in the next section. 
\bigskip

\begin{remark}
In \cite{Chapoton} and \cite{VanDerLaan}, F.~Chapoton and P.~van~der~Laan 
independently defined essentially the same group $G_{\QQ}$ of formal series 
for any algebraic operad $\QQ$ over $\Q$. We point out here the relationship 
between their construction and ours. 

The group $G_{\QQ}$ is formed of formal sums $\sum_{\mu\in\QQ} \mu$ 
with $\mu_1=\id$, endowed with the operation 
\begin{align}
\label{Q-composition}
\left(\sum_{\mu\in\QQ} \mu\right)\circ\left(\sum_{\eta\in\QQ} \eta\right)
&= \sum_{\mu\in\QQ}\ \ \sum_{\eta_1,\dots,\eta_{|\mu|}\in\QQ}
\gamma(\mu;\eta_1,\dots,\eta_{|\mu|}), 
\end{align}
where $\gamma$ denotes the operadic compositions in $\QQ$. 

In the case where $\QQ=\Q\P$ comes from a set-operad $\P$, the two 
constructions are essentially the same. 
Each vector space $\QQ(n)=\Q\P(n)$ has a canonical basis given by the 
elements of $\P(n)$, therefore any $\mu\in\QQ(n)$ can be written as a linear 
combination of these elements, that is $\mu=\sum_{p\in\P(n)} \mu_p\ p$, 
where $\mu_p\in\Q$ are scalar coefficients. 
Then the composition given by Eq.~(\ref{Q-composition}) is just the
linear extension of the composition given by Eq.~(\ref{P-composition}). 
In other words, the group $G_{\QQ}$ coincides with the group $\Gd_{\P}(\Q)$ 
of $\P$-expanded diffeomorphisms with coefficients in the ground field.  

The advantage of Chapoton-van~der~Laan's construction is that the group 
$G_{\QQ}$ can be defined for algebraic operads which are not set-operads. 
The advantage of our construction is that for set-operads we can distinguish 
between operadic elements and coefficients. These can then be chosen in any 
associative, commutative and unital algebra. 
\end{remark}

The construction of the group $\Gd_{\P}(A)$ is functorial in $\P$. 
The order map $|\ |:\P\longrightarrow\As$ is a canonical morphism of operads. 
A section is simply an operad morphism $\As\longrightarrow\P$. 
Let us call {\em associative\/} an element $p_2\in\P(2)$ such that 
$\gamma(p_2;p_2,\id)=\gamma(p_2;\id,p_2)$, and {\em multiplicative\/} 
an operad $\P$ equipped with an associative element. 
Then any associative $p_2\in\P(2)$ gives a section 
$n\mapsto p_n=\gamma(p_2;p_{n-1},\id)$. 
Therefore we can easily compare the series expanded over $\P$ with the usual 
ones. 

\begin{proposition}
\label{Gd_P->Gd}
For any set-operad $\P$, there is a canonical morphism of groups 
$\pi:\Gd_{\P}(A) \longrightarrow \Gd(A)$, induced by the order map $|\ |$. 
Moreover, if $\P$ is multiplicative, any associative $p_2\in\P(2)$ gives 
a section $i_{p_2}:\Gd(A) \longrightarrow \Gd_{\P}(A)$. 
\end{proposition}

Again, the construction of the group $\Gd_{\P}(A)$ is functorial also in $A$, 
and the group functor $\Gd_{\P}$ is represented by a commutative Hopf algebra 
$\Hd_{\P}$. As an algebra, $\Hd_{\P}$ is a polynomial algebra with generators 
indexed by $\PP$ and graded by the shifted degree. 
As shown by van~der~Laan for algebraic operads, the coordinate ring 
$\Hd_{\P}$ admits a non-commutative lift. 

Consider the free associative algebra $\Hdnc_{\P}=\Q\langle\PP\rangle$,  
graded by the shifted grading $\deg(u)=|u|-1$ for $u\in\P$, and where 
the element $\id\in\P(0)$ is identified with the formal unit $1$.  
Define a coproduct 
$\Dd_{\P}:\Hdnc_{\P}\longrightarrow\Hdnc_{\P}\otimes\Hdnc_{\P}$
as the algebra morphism given on the generators $u$ by 
\begin{align*}
\Dd_{\P}(u) &:= \underset{\gamma(p;q_1,\dots,q_{|p|})=u}
{\sum_{p,q_1,...,q_{|p|}\in\PP\cup\{\id\}}} 
p\otimes q_1\dots q_{|p|}. 
\end{align*}
Define also a counit $\varepsilon:\Hdnc_{\P}\longrightarrow\Q$ as the algebra 
morphism given on all $p\neq\id$ by $\varepsilon(p)=0$. 

\begin{proposition}
\label{Hd_P}
The algebra $\Hdnc_{\P}$ is a graded and connected Hopf algebra, 
neither commutative nor cocommutative. 

The coordinate ring of the group $\Gd_{\P}$ is the abelian 
quotient $\Hd_{\P}= \Q[\PP]$ of the Hopf algebra $\Hdnc_{\P}$. 
\end{proposition}

\begin{proof}
The fact that $\Hd_{\P}$ is the coordinate ring of the group $\Gd_{\P}(A)$ 
is obvious. The existence of a non-commutative lift is ensured by the 
assumption that the operad $\P$ is non-symmetric. In this case, in fact, 
the operadic composition fixes the order of the operations, and this 
guarantees the coassociativity of the coproduct lifted to tensor products.
\end{proof}
\bigskip 

\noindent{\bf Example.} 
The coordinate ring of the group functor $\Gd$ is the polynomial
algebra $\Q[a_1,a_2,\dots]$ on the graded generators $a_n$, one in
each degree $n \geq 1$. It is a Hopf algebra, with coproduct given by
\begin{align*} 
\Dd(a_n) &= \sum_{m=0}^n a_m \otimes \sum_{l=0}^{n-m} \binom{m+1}{l}
\underset{p_1+2p_2+\dots+(n-m)p_{n-m}=n-m}{\underset{p_1+p_2+\dots +p_{n-m}=l}
{\sum_{p_1,\dots,p_{n-m}\geq 0}}}
\frac{l!}{p_1!\dots p_{n-m}!}\ \ a_1^{p_1}\dots a_{n-m}^{p_{n-m}} ,
\end{align*}
counit $\varepsilon(a_n)=0$ for $n \geq 1$, and antipode defined recursively. 
Up to a rescaling of the generators, this Hopf algebra is known as the 
{\em Fa\`a di Bruno Hopf algebra\/}, cf. \cite{FaadiBruno}, \cite{JoniRota} 
or \cite{Majid}. 

The non-commutative version of this Hopf algebra is exactly the one 
defined in~\cite{BFK}, with coproduct lifted as 
\begin{align*} 
\Dd(a_n) &= \sum_{m=0}^n a_m \otimes 
\underset{k_0+k_1+\dots+k_m=n-m}{\sum_{k_0,k_1,\dots,k_m\geq 0}} 
a_{k_0} a_{k_1}\dots a_{k_m} 
\qquad (a_0=1). 
\end{align*}
\bigskip 

The non-commutative version of Proposition~\ref{Gd_P->Gd} gives the following 
result.  

\begin{proposition}
\label{Hd->Hd_P} 
There is a canonical morphism of Hopf algebras 
$\Hdnc\longrightarrow\Hdnc_{\P}$.
Moreover, if $\P$ is multiplicative, each associative $p_2\in\P(2)$ 
gives a section $\Hdnc_{\P}\longrightarrow\Hdnc$. 
\end{proposition}

\begin{proof}
The map $\Hdnc\longrightarrow\Hdnc_{\P}$ is explicitly given by 
\begin{align*}
a_n \longmapsto \sum_{p\in\P(n+1)} p. 
\end{align*} 
It is an inclusion if $\P(2)$ is not empty, and its section 
$\Hdnc_{\P}\longrightarrow\Hdnc$ is given by the map 
\begin{align*}
u\in \P(n) \longmapsto \left\{\begin{array}{ll} a_{n-1} & 
\mbox{if $u=p_n=\gamma(p_{n-1};p_2,\id,\dots,\id)$} \\ 
0 & \mbox{otherwise}\end{array}\right.. 
\end{align*}
The verification that these two maps are morphisms of Hopf algebras 
is trivial on the abelian quotients, and the passage to the non-commutative 
algebras is straightforward. 
\end{proof}


\subsection{Set-operad of trees and tree-expanded formal 
diffeomorphisms}
\label{tree-diffeomorphisms}

Any tree $t\in \Y$ can be written as a monomial in the vertex tree
$\treeA$, using the over and under products and suitable parentheses.
For instance,
\begin{align*}
\treeAB = \treeA\over\treeA, & \qquad \treeBA=\treeA\under\treeA, \\ 
\treeBAC = (\treeA\under\treeA)\over\treeA, & \qquad 
\treeACA = (\treeA\over\treeA)\under\treeA=\treeA\over(\treeA\under\treeA). 
\end{align*} 
This decomposition is not unique in general, as shown by the last example.
\bigskip 

For any tree $t\in \Y$, we call {\em $t$-product\/} the map 
$\mu_t:\Y^{\times |t|} \longrightarrow \Y$ which describes the
tree $t$ as an over and under product of $\treeA$ by itself\footnote{
  A similar monomial is considered by J.-L.~Loday in
  \cite{LodayATree}, based on the decomposition of a tree into some
  {\em left\/} and {\em right\/} products different from the over and
  under products considered here.}.  
More precisely, $\mu_t$ reproduces the tree $t$ when evaluated on $|t|$ 
copies of $\treeA$, that is $\mu_t(\treeA,\treeA,\dots,\treeA)=t$, 
and computes the product defined by the shape of $t$ on all the other trees
$s_1,s_2,\dots,s_{|t|} \in \Y$. Graphically, this means that, for any
trees $s_1,s_2,\dots,s_{|t|} \neq \treeO$, the tree
$\mu_t(s_1,s_2,\dots,s_{|t|})$ is obtained by replacing each internal
vertex of $t$, which has shape $\treeA$, by the tree $s_i$, in the
order given by the parentheses arising in the decomposition of $t$ by
$\treeA$.  The tree thus obtained clearly has order
\begin{align}
\label{order-monomial}
|\mu_t(s_1,s_2,\dots,s_{|t|})| & = |s_1|+|s_2|+\dots+|s_{|t|}|. 
\end{align}

In particular, if $t=\treeA$, the map $\mu_{\treeAsmall}:\Y \longrightarrow
\Y$ acts as the identity, that is $\mu_{\treeAsmall}(s)=s$ for any $s\neq
\treeO$. Other examples are:
\begin{align*}
&\treeBAC = (\treeA \under \treeA) \over \treeA  
\quad\mbox{hence}\quad 
\mu_{\treeBAC}(s_1,s_2,s_3) = (s_1 \under s_2) \over s_3 
= \lgraft{s_3}{\rgraft{s_1}{s_2}}, \\ 
&\treeACA = \treeA \over \treeA \under \treeA 
\quad\mbox{hence}\quad 
\mu_{\treeACA}(s_1,s_2,s_3) = s_1 \over s_2 \under s_3 
= \lgraft{\rgraft{s_2}{s_3}}{s_1} \quad . 
\end{align*}

We call {\em tree-product\/} the collection of the $t$-products 
given by 
\begin{align*}
\mu: \bigcup_{t\in \Y} \Y^{\times |t|} \longrightarrow \Y. 
\end{align*}

\begin{lemma}
\label{tree-operad}
The collection of trees $\Y=\cup_{n\geq 1} Y_n$ forms a set-operad, with 
operadic composition given by the tree-product $\mu$. 
\end{lemma}

\begin{proof}
Since $\Y=\cup_{m=1}^\infty Y_m$, we have 
\begin{align*}
\bigcup_{t\in \Y} \Y^{\times |t|} = 
\bigcup_{n=1}^\infty Y_n \times \Y^{\times n} = 
\underset{m_1,\dots,m_n \geq 1}{\bigcup_{n\geq 1}} 
Y_n \times Y_{m_1} \times\dots\times Y_{m_n}, 
\end{align*}
and we see in particular that for any choice $n\geq 1$ and 
$m_1,\dots,m_n\geq 1$, the map $\mu$ restricted to 
$Y_n \times Y_{m_1} \times\dots\times Y_{m_n}$ 
takes value in the homogeneous component $Y_{m_1+\dots+m_n}$ of $\Y$. 
In other words, the tree-product can be seen as the collection of the maps 
\begin{align*}
\mu_{n,m_1,\dots,m_n}: Y_n \times Y_{m_1} \times\dots\times Y_{m_n}
\longrightarrow Y_{m_1+\dots+m_n} 
\end{align*}
given on the trees $t\in Y_n$, $s_1\in Y_{m_1}$,\dots, $s_n\in Y_{m_n}$ by 
\begin{align*}
\mu_{n,m_1,\dots,m_n}(t;s_1,\dots,s_n) &= \mu_t(s_1,\dots,s_n). 
\end{align*}
The graphical interpretation of the tree-product given above ensures that 
the maps $\mu_{n,m_1,\dots,m_n}$ act by substitution of the variables 
(internal vertices) in $t$ by the operations $s_1,\dots,s_n$. 

Then it suffices to show that the tree-product $\mu$ is 
associative, in the sense that for any tree $t \in \Y$ with $|t|=n$, 
any choice of $n$ trees $s_1,\dots,s_n \in \Y$ with 
$|\mu_t(s_1,\dots,s_n)|=\sum_{i=1}^n |s_i|= m$, 
and any choice of $m$ trees $u_1,\dots,u_m \in \Y$, the two trees  
\begin{align}
\label{mu-asso1}
\mu_{\mu_t(s_1,\dots,s_n)}(u_1,\dots,u_m)
\end{align}
and
\begin{align}
\label{mu-asso2}
\mu_t \big(\mu_{s_1}(u_1,\dots,u_{|s_1|}), 
\mu_{s_2}(u_{|s_1|+1},\dots,u_{|s_1|+|s_2|}), \dots,
\mu_{s_n}(u_{|s_1|+\dots+|s_{n-1}|+1},\dots,u_m) \big)
\end{align}
coincide. 

To obtain the tree (\ref{mu-asso1}), we first construct the tree 
$\mu_t(s_1,\dots,s_n)$ by replacing each vertex of $t$ with each of the 
trees $s_1,\dots,s_n$, in the order given by $\mu_t$ as a monomial w.r.t. 
the over and under products. 
By Eq.~(\ref{order-monomial}), the tree thus obtained has exactly 
$\sum_{i=1}^n |s_i| =m$ vertices. Therefore, we can apply the 
$\mu_t(s_1,\dots,s_n)$-product to the $m$ trees $u_j$, and get the final 
tree $\mu_{\mu_t(s_1,\dots,s_n)}(u_1,\dots,u_m)$. 

Since the $\mu_t(s_1,\dots,s_n)$-product of $m$ trees contains all the 
sub-products of shapes $s_i$ delimited by parenthesis, and moreover they are 
all ordered by the shape of $t$, the final result is the same that we obtain 
if we first apply each $\mu_{s_i}$-product to the suitable package of trees 
$u_j$, and then apply the $t$-product to the $n$ new trees thus obtained. 
In summery, the resulting tree yields exactly (\ref{mu-asso2}). 
\end{proof}

The over and under operations on trees are associative operations which 
also satisfy one extra property: for any $s,t,u \in \Y$ (in fact $s$ and $u$ 
could be equal to $\treeO$), we have 
\begin{align*}
(s \over t)\under u = s \over (t\under u). 
\end{align*}
A set endowed with two associative operations veryfing this extra property 
was already considered by T.~Pirashvili in \cite{Pirashvili}, where it 
is called {\em duplex\/}. 
The operad which characterizes such operations was considered by M.~Aguiar 
and M.~Livernet in \cite{AguiarLivernet}, and by J.-L.~Loday in 
\cite{LodayTriple}, from whom we adopt the terminology. 

For our purpose, we then call {\em duplicial operad\/}, $\Dup$, 
the set-operad generated by the operations {\em over\/}, $\over$, and 
{\em  under\/}, $\under$. More precisely, $\Dup$ is the quadratic operad
obtained as the quotient of the free set-operad (with identity)
$\F=\F(\over,\under)$ on two binary operations, satisfying the three
relations
\begin{align*}
(a\over b)\over c &= a\over (b\over c)  \\ 
(a\over b)\under c &= a\over (b\under c) \\ 
(a\under b)\under c &= a\under (b\under c), 
\end{align*}
whenever the operations are applied to three elements $a,b,c$. 
The operadic composition on $\Dup$ is induced by that on $\F$, given, 
for any $n,m_1,\dots,m_n \geq 1$, by the map
\begin{align*}
\F(n)\times \F(m_1)\times\dots\times\F(m_n)
\longrightarrow \F(m_1+\dots +m_n)
\end{align*}
which sends the operations $(p,q_1,\dots,q_n)$ into the operation
obtained by inserting each operation $q_i$ into the $i$th variable
of the operation $p$.

\begin{theorem}
The set-operad of trees is isomorphic to the duplicial operad. 
\end{theorem}  

\begin{proof}
  Let us define a morphism $\kappa$ from the $\Dup$ operad to the operad 
  of trees, by sending the identity $\id\in\Dup(1)$ to the vertex tree 
  $\treeA$, and the generating operations $\over\mapsto \treeAB$ and
  $\under\mapsto\treeBA$. Since the trees $\treeAB$ and $\treeBA$ 
  satisfy the three necessary relations in the operad of trees, 
  the map $\kappa$ can be extended to a unique morphism of operads.
  To show that the morphism $\kappa$ is a bijection, it is sufficient 
  to prove that the free $\Dup$-monoid on one generator is isomorphic 
  to the set of trees, considered with the over and under products. 
  This was proved by J.-L.~Loday in \cite{LodayTriple} in the algebraic 
case. 
\end{proof}

Applying Theorem~\ref{Gd_P} to the duplicial operad, we obtain a group 
\begin{align*}
\Gd_Y(A) &:= \{ \varphi(x)=\sum_{t\in \Y} \varphi_t\ x^t, 
\ \varphi_t \in  A,\ \varphi_{\treeA}=1 \} 
\end{align*}
of tree-expanded formal diffeomorphisms. 
The composition of two tree-expanded series 
$\varphi(x)=\sum_{t} \varphi_t\ x^t$ and $\psi(x)=\sum_{s} \psi_s\ x^s$ 
is given by 
\begin{align*}
(\varphi \circ \psi)(x) 
& = \underset{s_1,s_2,\dots,s_{|t|}\in \Y}{\sum_{t\in \Y}} 
\varphi_t\ \psi_{s_1} \psi_{s_2}\dots \psi_{s_{|t|}} \ 
x^{\mu_t(s_1,s_2,\dots,s_{|t|})}. 
\end{align*}
If we define the power of the series $\psi(x)$ by a tree $t\neq\treeO$ as  
$\psi(x)^t = \mu_t(\psi(x),\psi(x),\dots,\psi(x))$, 
the composition of tree-expanded series can also be seen as 
a substitution, that is $(\varphi \circ \psi)(x) = \varphi(\psi(x)) 
= \sum_{t\in \Y} \varphi_t\ \psi(x)^t$. 
\bigskip 

\noindent{\bf Example.} 
Let $\varphi(x) = x^{\treeA} + a\ x^{\treeAB} + b\ x^{\treeBA}$ and 
$\psi(x)=x^{\treeA}+c\ x^{\treeAB}+d\ x^{\treeBA}$, with $a,b,c,d \in  A$. 
Since $\mu_{\treeAB}(t,s)=t\over s$ and $\mu_{\treeBA}(t,s)=t\under s$, 
we compute  
\begin{align*}
(\varphi \circ \psi)(x) 
&= \psi(x) + a\ \psi(x)^{\treeAB} + b\ \psi(x)^{\treeBA} \\ 
&= \psi(x) + a\ \psi(x)\over \psi(x) + b\ \psi(x)\under \psi(x) \\ 
&= x^{\treeA} + (a+c)\ x^{\treeAB} + (b+d)\ x^{\treeBA} \\ 
&\qquad + 2ac\ x^{\treeABC} + ad\ x^{\treeBAC} + (ad+bc)\ x^{\treeACA} 
+ bc\ x^{\treeCAB} + 2bd\ x^{\treeCBA} \\ 
&\qquad + ac^2\ x^{\treeABCD} + acd\ x^{\treeBACD} + acd\ x^{\treeABDA} 
+ ad^2\ x^{\treeBADA} + bc^2\ x^{\treeADAB} \\ 
&\qquad + bcd\ x^{\treeADBA} + bcd\ x^{\treeDCAB} + bd^2\ x^{\treeDCBA}. 
\end{align*}

More interesting examples of compositions of tree-expanded diffeormorphisms 
are computed by F.~Chapoton in \cite{ChapotonFrabetti}. 
\bigskip

Since $\Dup(2)$ is not empty, applying Proposition~\ref{Gd_P->Gd}, we see
that the order map $\pi$ of Section~\ref{tree-invertible} gives also a
surjective morphism of groups $\pi:\Gd_Y(A)\longrightarrow\Gd(A)$, sending
$x^t$ to $x^{|t|}$.

Vice versa, since $\Dup(2)=\{\over,\under\}$ contains two operations
which are both associative, the projection $\pi$ has two sections, the
maps $\io,\iu:\Gd(A)\longrightarrow\Gd_Y(A)$ sending $x$ to $x^{c_n}$, where
$c_n$ denotes respectively the left and the right comb trees of order
$n$. Therefore the two subgroups of $\Gd_Y(A)$ containing the
tree-expanded diffeomorphisms expanded only on the comb trees are
isomorphic to the group of usual formal diffeomorphisms.  

In summery, the relationship between formal diffeomorphisms expanded 
over trees and over natural numbers can be deduced by the natural maps 
between the associated operads: 
\begin{equation*}
\label{diagramme_operades}
\xymatrix{
\As \ar @<.7ex>[r]^{\io}\ar @<-.7ex>[r]_{\iu} &\Dup \ar[r] 
& \Dias \ar[r]^{\pi} & \As}.
\end{equation*}
\bigskip

The aim of introducing tree-expanded formal diffeomorphisms is to describe 
the group law of the group $\Ga(A)$ dual to the ``charge renormalization 
Hopf algebra'' $\Ha$ expressed by means of tree-expanded series. 
The group $\Gd_Y(A)$ indeed allows to generalize the composition of usual 
series to the tree-expanded series, but its coordinate ring is much bigger 
than the Hopf algebra $\Ha$. 
In fact, if we denote by $\YY=\cup_{n\geq 2} Y_n$ the set of all trees but 
$\treeO$ and $\treeA$, the coordinate ring of the group functor $\Gd_Y$ 
is given by the polynomial algebra $\Hd_Y=\Q[\YY]\cong\Q[\Y]/(\treeA-1)$ 
generated by all the trees of order at least $2$, instead of only a subset 
of them. 
The interest of the algebra $\Hd_Y$ is that it allows a non-commutative lift. 

\begin{corollary}
The free associative algebra $\Hdnc_Y=\Q\langle\Y\rangle/(\treeA-1)$ 
is a graded and connected Hopf algebra, with respect to the shifted grading 
$\deg(t)=|t|-1$. The coproduct is defined on the generators $u\in \Y$ as
\begin{align*}
\Dd_Y(u)
& :=\underset{u=\mu_t(s_1,\dots,s_{|t|})}{\underset{s_1,\dots,s_{|t|}\in \Y}
{\sum_{t\in \Y}}} t \otimes s_1\dots s_{|t|}, 
\end{align*}
and the counit is $\varepsilon(u)=0$ for any $u\in \Y$, $u\neq \treeA$. 
This Hopf algebra is neither commutative nor cocommutative. 
\end{corollary}

Note in particular that $\Dd_Y(\treeA) =\treeA \otimes \treeA$, because 
$\treeA = \mu_{\treeA}(\treeA)$, therefore $\treeA$ is a group-like element 
and can be identified with the unit. 
Setting $\treeA=1$, the coproduct on small trees is: 
\begin{align*}
\Dd_Y(\treeAB) &= \treeAB \otimes 1 + 1 \otimes \treeAB, \\
\Dd_Y(\treeBA) &= \treeBA \otimes 1 + 1 \otimes \treeBA, \\
\Dd_Y(\treeABC) &= \treeABC \otimes 1 + 2 \treeAB \otimes \treeAB +
1 \otimes \treeABC, \\
\Dd_Y(\treeBAC) &= \treeBAC \otimes 1 + \treeAB \otimes \treeBA
+ 1 \otimes \treeBAC, \\
\Dd_Y(\treeACA) &= \treeACA \otimes 1 + \treeAB \otimes\treeBA 
+ \treeBA \otimes \treeAB + 1 \otimes \treeACA, \\
\Dd_Y(\treeCAB) &= \treeCAB \otimes 1 + \treeBA \otimes \treeAB 
+ 1 \otimes \treeCAB, \\
\Dd_Y(\treeCBA) &= \treeCBA \otimes 1 + 2 \treeBA \otimes
\treeBA + 1 \otimes \treeCBA. 
\end{align*} 

By Proposition~\ref{Hd->Hd_P}, the map 
\begin{align*}
\Hdnc \longrightarrow \Hdnc_Y: a_n \longmapsto \sum_{|t|=n+1} t
\end{align*}
gives an inclusion of Hopf algebras. 


\section{Action of tree-expanded diffeomorphisms on 
tree-expanded invertible series}
\label{section-action} 

In this section we describe an action of the group of tree-expanded 
diffeomorphism on that of invertible series which will be used in the 
next section to construct the group $\Ga$.  

Such an action can be defined on the analogue groups of series expanded 
over any set-operad which admits a canonical associated graded monoid. 
We begin by discussing the general case, and then specify to the duplicial 
operad. 

\subsection{Groups of series expanded over an associative 
set-operad}

From now on, we assume that $\P$ is a multiplicative set-operad, and 
we denote by $p_2$ the associative element in $\P(2)$. 
This is equivalent to require that there exists an operad morphism 
$\As\longrightarrow\P$. 

Then we can naturally associate to $\P$ a graded monoid by setting 
$\M_{\P}=\P\cup\{e\}$, with $\M_{\P}(n)=\P(n)$ for $n\geq 1$ and 
$\M_{\P}(0)=\{e\}$, where $e$ is a formal element which is taken as the 
neutral element. The multiplication of $p\in\M_{\P}(m)$ and $q\in\M_{\P}(n)$ 
is defined by 
\begin{align*}
p\cdot q &:= \gamma_{\P}(p_2;p,q) \in\M_{\P}(m+n). 
\end{align*}

The graded monoid $\M_{\P}$ then determines the algebra $A[[\M_{\P}]]$
and the group $\Gi_{\M_{\P}}(A)$. 
The composition of $\P$-expanded series given in Eq.~\ref{P-composition} 
can be extended to a map 
$A[[\M_{\P}]]\times A[[\P]]\longrightarrow A[[\M_{\P}]]$, by setting
\begin{align*}
(f,\varphi)&\longmapsto f^{\varphi}(x)
:= f_e\ x^e + \sum_{p\in\P}\ \ \sum_{q_1,\dots,q_{|p|}\in\P}
f_{p}\ \varphi_{q_1}\dots \varphi_{q_n}\ x^{\gamma(p;q_1,\dots,q_{|p|})}. 
\end{align*}
Restricting this map on one side to the group $\Gi_{\M_{\P}}(A)$ of 
$\M_{\P}$-expanded invertible series and on the other side to the group 
$\Gd_{\P}(A)$ of $\P$-expanded diffeomorphisms, we obtain the following 
result. 

\begin{theorem}
\label{P-action-series}
The composition defines a graded right action
$\Gi_{\M_{\P}}(A)\times\Gd_{\P}(A)\longrightarrow\Gi_{\M_{\P}}(A)$ 
of the group of $\P$-expanded diffeomorphisms on the group of 
$\M_{\P}$-expanded invertible series.
\end{theorem}

\begin{proof}
The compatibility of the action with the composition in $A[[\P]]$, that is
$(f^{\varphi})^{\psi} = f^{\varphi\circ\psi}$, is ensured by the 
associativity of the operadic composition $\gamma$. 
The computations are exactly the same as those which show that the 
composition $\circ$ in $A[[\P]]$ is associative, cf. Theorem~\ref{Gd_P}.

We show that the action preserves the multiplication $\cdot$ in 
$A[[\M_{\P}]]$, that is $(f \cdot g)^{\varphi}=f^{\varphi}\cdot g^{\varphi}$. 
Let us compute the two terms of this equality separately.
Since the action is additive, on one side we have
\begin{align*}
(f \cdot g)^{\varphi}(x) &= f_e\ g_e\ x^e + g_e\ f^\varphi(x)+ f_e\ g^\varphi
+\underset{u_1,\dots,u_{|p\cdot q|}\in\P}{\sum_{p,q\in\P}} 
f_p\ g_q\ \varphi_{u_1}\dots\varphi_{u_{|p \cdot q|}}\ 
x^{\gamma(p\cdot q;u_1,\dots,u_{|p\cdot q|})},
\end{align*}
and on the other side
\begin{align*}
(f^{\varphi}\cdot g^{\varphi})(x) &= f_e\ g_e\ x^e + g_e\, f^\varphi(x) 
+ f_e\, g^\varphi(x) \\ 
&\quad 
+\underset{u_{|p|+1},\dots,u_{|p|+|q|}\in\P}{\underset{u_1,\dots,u_{|p|}\in\P}
{\sum_{p,q\in\P}}} f_p\ g_q\ \varphi_{u_1}\dots\varphi_{u_{|p|+|q|}}\
x^{\gamma(p;u_1,\dots,u_{|p|})\cdot \gamma(q;u_{|p|+1},\dots,u_{|p|+|q|})}.
\end{align*}
Since $|p\cdot q|=|p|+|q|$, the two terms coincide if the exponents
coincide, and this is again ensured by the associativity of the
operadic composition. 

Finally, the action is graded with respect to the two different gradings 
given on $A[[\M_{\P}]]$ and $A[[\P]]$, namely $\deg(p)=|p|$ if 
$p$ is in the monoid $\M_{\P}$ and $\deg(p)=|p|-1$ if $p$ is in the 
set-operad $\P$. In fact, if $p\in\M_{\P}$ is different from $e$, then for any 
$q_1,\dots,q_{|p|}\in\P$, we have
\begin{align*}
\deg(p) + \deg(q_1 \dots q_{|p|}) 
&= |p| + \sum_{i=1}^{|p|} (|q_i|-1) = \sum_{i=1}^{|p|} |q_i| \\ 
&= \deg(\gamma(p;q_1,\dots,q_{|p|})).   
\end{align*}
\end{proof}

\noindent{\bf Example.} 
In the case $\P=\As$, we have $\M_{\As}=\N$ and the right action 
$\Gi(A)\times\Gd(A)\longrightarrow\Gi(A)$ is the ordinary composition 
$(f^{\varphi})(x)= f(\varphi(x))$ of invertible series by formal 
diffeomorphisms. 
\bigskip 

Given a right action of $\Gd_{\P}(A)$ on the group $\Gi_{\M_{\P}}(A)$, 
we can define the semi-direct product $\Gd_{\P}(A) \ltimes \Gi_{\M_{\P}}(A)$ 
in the ususal way: as a set we take $\Gd_{\P}(A) \times \Gi_{\M_{\P}}(A)$, 
and the group law is given by
\begin{align*}
(\varphi,f) \cdot (\psi,g) &= (\varphi\circ\psi,f^{\psi}\cdot g), 
\end{align*}
for any $\varphi,\psi\in\Gd_{\P}(A)$ and $f,g\in\Gi_{\M_{\P}}(A)$. 
The order map $\pi:A[[\P]]\longrightarrow A[[x]]$ sending 
$x^p$ to $x^{|p|}$ then induces a projection of groups 
$\Gd_{\P}(A) \ltimes \Gi_{\M_{\P}}(A) \longrightarrow \Gd(A)\ltimes\Gi(A)$. 
\bigskip

The right action of $\Gd_{\P}(A)$ on $\Gi_{\M_{\P}}(A)$ becomes a right 
coaction on the coordinate rings. It can be lifted to the non-commutative Hopf 
algebras, as follows. Let 
$\dinv_{\P}: \Hinc_{\M_{\P}}\longrightarrow \Hinc_{\M_{\P}}\otimes\Hdnc_{\P}$ 
be the algebra morphism given on the generators $u\in\P$ by 
\begin{align*}
\dinv_{\P}(u) &:= \underset{\gamma(p;q_1,\dots,q_{|p|})=u}
{\underset{q_1,...,q_{|p|}\in\PP\cup\{\id\}}{\sum_{p\in\P}}} 
p \otimes q_1\dots q_{|p|},
\end{align*}
and on the unit by $\dinv_{\P}(e) = e \otimes \id$. 
Note that formally $\dinv_{\P}(u)=\Dd_{\P}(u)$ on the generators $u\in\P$, 
but these two elements have different degrees in their proper spaces, 
as well as $u$ itself. 

\begin{proposition}
\label{dinv_P}
The map $\dinv_{\P}$  is a graded right coaction of the Hopf algebra 
$\Hdnc_{\P}$ on the algebra $\Hinc_{\M_{\P}}$. 
Moreover, $\Hinc_{\M_{\P}}$ is an $\Hdnc_{\P}$-comodule coalgebra. 
\end{proposition}

\noindent
In other words, the map $\dinv_{\P}$ satisfies the two following conditions: 
\begin{align*}
(\dinv_{\P}\otimes\Id)\ \dinv_{\P} 
&= (\Id\otimes\Dd_{\P})\ \dinv_{\P}, \\ 
(\Di_{\M_{\P}} \otimes\Id)\ \dinv_{\P} 
&= (\Id\otimes\Id\otimes m)\ (\Id\otimes\tau\otimes\Id)\ 
(\dinv_{\P}\otimes\dinv_{\P})\ \Di_{\M_{\P}}, 
\end{align*}
where $m$ denotes the multiplication in the algebra $\Hdnc_{\P}$, and
$\tau$ denotes the twist.
\bigskip 

\noindent{\bf Example.} 
The case $\P=\As$ gives a right coaction 
$\dinv:\Hi\longrightarrow\Hi\otimes\Hd$ which sends a generator $b_n$ into 
\begin{align*}
\dinv(b_n) &= \sum_{m=0}^n b_m \otimes 
\underset{k_1+\dots+k_m=n-m}{\sum_{k_1,\dots,k_m\geq 0}} 
a_{k_1}\dots a_{k_m}, 
\qquad (a_0=1). 
\end{align*}
The non-commutative version 
$\dinv:\Hinc\longrightarrow\Hinc\otimes\Hdnc$ was considered in \cite{BFK},
where it was denoted by $\ddif$. In this paper we reserve the
notation $\ddif$ for a coaction of $\Hd$ on itself, which will be
introduced in Section~\ref{section-subgroup}. 
\bigskip

The group functor $\Gd_{\P}\ltimes\Gi_{\M_{\P}}$ is represented 
by the semi-direct coproduct (or smash coproduct) Hopf algebra 
$\Hd_{\P}\ltimes\Hi_{\M_{\P}}$. 
As an algebra, this is the tensor product $\Hd_{\P}\otimes\Hi_{\M_{\P}}$. 
As a coalgebra, it is endowed with the twisted coproduct defined 
on the generators $p\in\Hd_{\P}$ and $q\in\Hi_{\M_{\P}}$ by 
\begin{align*}
\Dl_{\P}(p\otimes q) 
&= \Dd_{\P}(p)\cdot [(\dinv_{\P}\otimes\Id)\ \Di_{\M_{\P}}(q)]. 
\end{align*}

Applying the results found by R.~Molnar in~\cite{Molnar}, we know that this 
Hopf algebra admits a non-commutative lift given by the semi-direct coproduct 
Hopf algebra $\Hd_{\P}\ltimes\Hinc_{\M_{\P}}$.

Instead, note that the semi-direct coproduct $\Hdnc_{\P}\ltimes\Hinc_{\M_{\P}}$
is at the same time an algebra and a coalgebra, but not a Hopf algebra 
because the non-commutativity of the algbera $\Hdnc_{\P}$ prevents 
the coproduct $\Dl_{\P}$ to be an algebra morphism.


\subsection{Tree-expanded series and actions}
\label{subsection-action-tree}

In this section we briefly illustrate the above results on the example 
of trees, using the duplicial operad.
\bigskip 

The duplicial operad has two binary operations, $\over$ and $\under$,
both associative. The graded monoid $\M_{\Dup}$ associated to these
operations are exactly the monoids of trees introduced in
Section~\ref{tree-invertible}. Therefore the group $\Gd_Y(A)$ acts on the
groups $\Go(A)$ and $\Gu(A)$, by composition. 

The action is compatible with the group structures of $\Go(A)$ and $\Gu(A)$,
therefore the semi-direct products $\Gd_Y(A) \ltimes \Go(A)$ and $\Gd_Y(A)
\ltimes \Gu(A)$ form two groups.  
\bigskip 

Proposition~\ref{dinv_P} then tells us how this action is reflected on the 
dual Hopf algebras on trees. 
Let us denote by $\Hi_Y = \Q[\Y] \cong \Q[Y](\treeO-1)$ the
coordinate ring of the proalgebraic set $\Gi_Y(A)$ of tree-expanded invertible
series, as introduced in Section~\ref{tree-invertible}, and by 
by $\Hinc_Y = \Q\langle \Y \rangle$ its non-commutative lift. 
Then $\Hp$ and $\He$ are the two Hopf algebras with underlying algebra 
$\Hinc_Y$ endowed with the ``pruning coproducts'' $\Dip$ and $\Die$.

\begin{corollary}
The algebra homomorphism $\dinv_Y:\Hinc_Y\longrightarrow\Hinc_Y\otimes\Hdnc_Y$ 
defined on the generators $u \in \Y$ formally as the coproduct $\Dd_Y$, that is
\begin{align*}
\dinv_Y(u)
&:=\underset{u=\mu_t(s_1,\dots,s_{|t|})}{\underset{s_1,\dots,s_{|t|}\in \Y}
{\sum_{t\in \Y}}} t \otimes s_1\dots s_{|t|}. 
\end{align*}
and which respects the units, that is $\dinv_Y(\treeO)=\treeO\otimes\treeA$, 
is a graded right coaction of the Hopf algebra $\Hdnc_Y$ on the algebra 
$\Hinc_Y$. 
\end{corollary}

The map induced by $\dinv_Y$ on the abelian quotients $\Hd_Y$ and $\Hi_Y$ 
is dual to the right action of the group $\Gd_Y(A)$ on the set $\Gi_Y(A)$.
\bigskip 

Note that the degrees and the units are different in the two algebras 
$\Hinc_Y$ and $\Hdnc_Y$. Therefore, even if the coproduct $\Dd_Y$ and 
the coaction $\dinv_Y$ are formally defined by the same expression, 
the meaning of the result is different. 
To see this difference, compare the value of $\Dd_Y$ on small trees, 
given at the end of Section~\ref{tree-diffeomorphisms}, with the 
following values of $\dinv_Y$, obtained by setting $\treeA=1$ in $\Hd_Y$, 
and $\treeO=1$ in $\Hi_Y$: 
\begin{align*}
\dinv_Y(\treeA) &= \treeA \otimes 1, \\ 
\dinv_Y(\treeAB) &= \treeAB \otimes 1 + \treeA \otimes \treeAB, \\
\dinv_Y(\treeBA) &= \treeBA \otimes 1 + \treeA \otimes \treeBA, \\
\dinv_Y(\treeABC) &= \treeABC \otimes 1 + 2 \treeAB \otimes \treeAB +
\treeA \otimes \treeABC, \\
\dinv_Y(\treeBAC) &= \treeBAC \otimes 1 + \treeAB \otimes \treeBA
+ \treeA \otimes \treeBAC, \\
\dinv_Y(\treeACA) &= \treeACA \otimes 1 + \treeAB \otimes\treeBA 
+ \treeBA \otimes \treeAB + \treeA \otimes \treeACA, \\
\dinv_Y(\treeCAB) &= \treeCAB \otimes 1 + \treeBA \otimes
\treeAB + \treeA \otimes \treeCAB, \\
\dinv_Y(\treeCBA) &= \treeCBA \otimes 1 + 2 \treeBA \otimes \treeBA 
+ \treeA \otimes \treeCBA. 
\end{align*} 

Proposition~\ref{dinv_P} tells us that the coaction $\dinv_Y$ 
is compatible with the coproducts $\Dip$ and $\Die$, and with the counit 
$\varepsilon$. Therefore $\Hp$ and $\He$ are coalgebra comodules over 
$\Hdnc_Y$. 
Then, the semi-direct coproduct algebras $\Hd_Y \ltimes \Hp$ and
$\Hd_Y \ltimes \He$ are non-commutative Hopf algebras, which lift the 
coordinate rings of the group functors $\Gd_Y\ltimes\Go$ and 
$\Gd_Y\ltimes\Gu$ respectively. 

Finally, the maps $b_n \longmapsto \sum_{|t|=n} t$ and 
$a_n\longmapsto \sum_{|t|=n+1} t$ define an inclusion of the Hopf algebra 
$\Hd\ltimes\Hinc$ into respectively $\Hd_Y\ltimes\Hp$ and $\Hd_Y\ltimes\He$. 


\section{Subgroup dual to the QED charge Hopf 
algebra}
\label{section-subgroup}

The renormalization of the electric charge in quantum electrodynamics 
was described in \cite{BFqedtree} by a commutative Hopf algebra $\Ha$ on trees
which was proved in \cite{BFK} to contain the Fa\`a di Bruno Hopf algebra, 
that is, the coordinate ring of $\Gd$, and which is different from $\Hd_Y$. 
Since this latter is the largest Hopf algebra on trees describing the 
composition of tree-expanded series, it is natural to look for a 
subgroup of $\Gd_Y$ having $\Ha$ as coordinate ring. 
We describe it in this section. 

To do it, we first introduce some intermediate subgroups of tree-expanded 
diffeomorphisms which exist for any multiplicative set-operad.
The final construction of the group $\Ga$, dual to $\Ha$, is possible only 
if the set-operad has two distinct associative elements with some suitable 
compatibility relation. 
At this level we specify the construction to the duplicial set-operad. 


\subsection{Subgroups of series expanded over a graded monoid 
set-operad}

Let $\P$ be a multiplicative set-operad and let $\M_{\P}$ be its associated 
graded monoid. 
For any associative, commutative and unital algebra $A$, we consider 
the two linear maps 
$\lambda,\rho:A[[\M_{\P}]]\longrightarrow A[[\P]]$ defined 
on a series $f(x)=f_e\ x^e+\sum_{p\in\P} f_p\ x^p$ by 
\begin{align*}
\lambda_f(x)&:=x^{\id}\cdot f(x) 
= f_e\ x^{\id}+\sum_{p\in\P} f_p\ x^{\id\cdot p}\\
\rho_f(x)&:=f(x) \cdot x^{\id} 
= f_e\ x^{\id}+\sum_{p\in\P} f_p\ x^{p\cdot\id}. 
\end{align*}
These maps are injective, and we denote their images in $A[[\P]]$ 
by $x^{\id}\cdot A[[\M_{\P}]]$ and by $A[[\M_{\P}]]\cdot x^{\id}$ 
respectively. 

\begin{theorem}
\label{subgroups}
The two sets
\begin{align*}
\Gl_{\P}(A):=x^{\id}\cdot\Gi_{\M_{\P}}(A) \qquad\mbox{and}\qquad 
\Gr_{\P}(A):=\Gi_{\M_{\P}}(A)\cdot x^{\id}
\end{align*}
are subgroups of the group $\Gd_{\P}(A)$. 
\end{theorem}

\begin{proof}
It suffices to show that the images of $\lambda$ and $\rho$ are stable 
under the composition of series. 
Given $\lambda_f(x)=x^{\id}\cdot f(x)$ and $\lambda_g(x)=x^{\id}\cdot g(x)$ 
in $x^{\id}\cdot A[[\M_{\P}]]$, we have to show that 
there exists an $h\in A[[\M_{\P}]]$ such that
$\left(\lambda_f \circ \lambda_g\right)(x)= x^{\id}\cdot h(x)$.
Using the compatibility of the composition with the multiplication
proved in Theorem~\ref{P-action-series}, we have
\begin{align*}
\left(\lambda_f \circ \lambda_g\right)(x) 
&= \left(x^{\id} \cdot f(x)\right)^{\lambda_g(x)} 
= \lambda_g(x) \cdot f^{\lambda_g}(x) 
= [x^{\id} \cdot g(x)] \cdot f^{\lambda_g}(x) \\
& = x^{\id} \cdot \left[g \cdot f^{\lambda_g} \right](x). 
\end{align*}
Therefore $\lambda_f\circ\lambda_g=\lambda_h$ if we set 
$h=g\cdot f^{\lambda_g}$.
Similarly, if $\rho_f$ and $\rho_g$ belong to $A[[\M_{\P}]]\cdot x^{\id}$, 
we have $\rho_f \circ \rho_g=\rho_h$, with $h=f^{\rho_g} \cdot g$.
\end{proof}
\bigskip 

The map $\lambda:\Gi_{\M_{\P}}(A)\longrightarrow\Gl_{\P}(A)$ is an 
isomorphism of sets and its inverse 
$\lambda^{-1}:\Gl_{\P}(A)\longrightarrow\Gi_{\M_{\P}}(A)$ sends a 
$\P$-expanded diffeomorphism of the form $\lambda_f(x)=x^{\id}\cdot f(x)$ 
to the $\M_{\P}$-expanded invertible series $f(x)$. 
Note that $\lambda$ is {\em not\/} a morphism of groups. 
Instead, its inverse $\lambda^{-1}$ is a 1-cocycle of $\Gl_{\P}(A)$ 
with values in $\Gi_{\M_{\P}}(A)$, with respect to the right action by 
composition, that is
\begin{align*} 
\lambda^{-1}(\psi) [\lambda^{-1}(\varphi\circ\psi)]^{-1} 
\lambda^{-1}(\varphi)^{\psi} &=x^e,
\end{align*} 
for any $\varphi,\psi\in\Gl_{\P}(A)$. 
\bigskip 

\noindent{\bf Example.} 
In the case $\P=\As$, the multiplication in the graded monoid 
$\M_{\As}=\N$ is commutative, therefore $\lambda=\rho$ and 
$\Gl_{\As}(A)=\Gr_{\As}(A)$. 
Moreover, the map $\lambda$ simply brings an invertible series $f(x)$ into 
$\lambda_f(x)=xf(x)$. This map is invertible on the whole space 
$A[[\As]]=xA[[x]]$, and its inverse $\lambda^{-1}$ brings a formal 
diffeomorphism $\varphi$ into the invertible series $\frac{\varphi(x)}{x}$. 
Therefore the group $\Gl_{\As}(A)$ coincides with the whole group $\Gd(A)$ 
of formal diffeomorphisms. 
\bigskip

In general, the two groups $\Gl_{\P}(A)$ and $\Gr_{\P}(A)$ are not
isomorphic, because the multiplication by $\id$ in the monoid 
$\M_{\P}$ is not commutative in general. 
\bigskip 

Consider now the map 
$(x^{\id}\cdot A[[M_{\P}]])\times A[[\P]]\longrightarrow 
x^{\id}\cdot A[[M_{\P}]]$ defined by 
\begin{align}
\label{P-lambda-action}
{\lambda_f}^{\psi}(x) &:= \lambda_{f^{\psi}}(x) = x^{\id}\cdot f^{\psi}(x) 
\nonumber \\ 
&= f_e\ x^{\id} + \underset{q_1,\dots,q_{|p|}\in\P}{\sum_{p\in\P}}
f_{p}\ \psi_{q_1}\dots\psi_{q_n}\ 
x^{\id\cdot\gamma(p;q_1,\dots,q_{|p|})}, 
\end{align}
for $f(x)=f_e\ x^e+\sum_{p\in\P} f_p\ x^p$ and 
$\psi(x)=\sum_{q\in\P} \psi_q\ x^q$. 
Similarly, define a map 
$(A[[M_{\P}]]\cdot x^{\id})\times A[[\P]]\longrightarrow 
A[[M_{\P}]]\cdot x^{\id}$, by setting 
\begin{align*}
{\rho_f}^\psi(x) &:= \rho_{f^{\psi}}(x)= f^{\psi}(x)\cdot x^{\id}. 
\end{align*}

\begin{theorem}
\label{Gd-action-Glr}
The map $(\lambda_f,\psi)\longmapsto {\lambda_f}^\psi= \lambda_{f^{\psi}}$,  
restricted to $\Gl_{\P}(A)\times\Gd_{\P}(A)\longrightarrow\Gl_{\P}(A)$,  
is a right action of $\Gd_{\P}(A)$ on $\Gl_{\P}(A)$. 
The analogue statement holds for $\Gr_{\P}(A)$. 
\end{theorem}

\begin{proof}
In fact the series $x^{\id}\in\Gd_{\P}(A)$ obviously acts as the identity, 
and for any $\psi,\eta\in\Gd_{\P}(A)$ we have 
\begin{align*}
({\lambda_f}^\psi)^\eta 
&= (\lambda_{f^{\psi}})^{\eta} = \lambda_{(f^{\psi})^{\eta}} 
= \lambda_{f^{\psi\circ\eta}} = {\lambda_f}^{\psi\circ\eta}. 
\end{align*}
\end{proof}

Note that the right action of $\Gd_{\P}(A)$ on $\Gl_{\P}(A)$ is indeed 
{\em different\/} from the composition of $\Gl_{\P}(A)$ by $\Gd_{\P}(A)$, 
where $\Gl_{\P}(A)$ is seen as a subgroup of $\Gd_{\P}(A)$. 
In other words, for any $\lambda_f\in \Gl_{\P}(A)$, and any 
$\psi\in\Gd_{\P}(A)$, in general we have 
${\lambda_f}^{\psi} \neq \lambda_f\circ\psi$. 
This is easily seen if we restrict the action of $\Gd_{\P}(A)$ to its subgroup 
$\Gl_{\P}(A)$. In this case, in fact, we have 
\begin{align*}
{\lambda_f}^{\lambda_g}=\lambda_{f^{\lambda_g}} 
&\neq \lambda_{g\cdot f^{\lambda_g}} = \lambda_f\circ\lambda_g. 
\end{align*}

\noindent{\bf Example.} 
In the case $\P=\As$, the action of $\Gd(A)$ on itself induced by 
the action of $\Gd(A)$ on $\Gi(A)$ has the following explicit form, 
\begin{align}
\label{dif-action}
\varphi^{\psi}(x) &= x+ \sum_{n=2}^\infty 
\left(\sum_{m=2}^n\ \ \underset{k_2,\dots,k_m\geq 1}{\sum_{k_2+\dots+k_m=n-1}} 
\varphi_m\ \psi_{k_2}\dots \psi_{k_m} \right)\ x^n, 
\end{align}
for any $\varphi,\psi\in\Gd(A)$. Comparing this expression with that 
of the composition given in Eq.~(\ref{composition}), we see that 
we obtain (\ref{dif-action}) if we set $\psi_{k_1}=\psi_1=1$ 
in (\ref{composition}). 
\bigskip

As usual, the construction of the groups $\Gl_{\P}(A)$ and $\Gr_{\P}(A)$ 
is functorial in $A$. Again, the coordinate rings $\Hl_{\P}$ and $\Hr_{\P}$ 
of the group functors $\Gl_{\P}$ and $\Gr_{\P}$ admit a straightforward 
non-commutative lift, as well as the actions by $\Gd_{\P}$. 
In Section~\ref{subsection-charge} we will use in particular the Hopf algebra 
$\Hr_{\P}$. Let us then describe explicitly only its non-commutative lift. 

Consider the free associative algebra 
$\Hrnc_{\P}=\Q\langle\M_{\P}\rangle/(e-1)\cong\Q\langle\P\rangle$, 
with grading given by the order, $\deg(u)=|u|$ for $u\in\P$, and where the 
element $e\in\M_{\P}(0)$ is identified with the unit $1$. 
Define a coproduct 
$\Dr_{\P}:\Hrnc_{\P}\longrightarrow\Hrnc_{\P}\otimes\Hrnc_{\P}$ 
as the algebra morphism given on the generators $u$ by 
\begin{align}
\label{Dr_P}
\Dr_{\P}(u) &:= 1\otimes u + 
\underset{u=\gamma(p;q_1\cdot\id,\dots,q_{|p|}\cdot\id)\cdot q_{|p|+1}}
{\underset{q_1,...,q_{|p|+1}\in\M_{\P}}{\sum_{p\in\P}}} 
p\otimes q_1\cdots q_{|p|+1}. 
\end{align}
Define a counit $\varepsilon:\Hrnc_{\P}\longrightarrow\Q$ as the algebra 
morphism given on any $u\in\P$ by $\varepsilon(u)=0$. 

\begin{theorem}
\label{Hr_P}
The algebra $\Hrnc_{\P}$ is a graded and connected Hopf algebra, neither 
commutative nor cocommutative. 

The coordinate ring of the group $\Gr_{\P}$ is the abelian 
quotient $\Hr_{\P}= \Q[\P]$ of the Hopf algebra $\Hrnc_{\P}$. 
\end{theorem}

\begin{proof}
The only difficulty of this result is the explicit form (\ref{Dr_P}) 
of the coproduct. 

In fact, if $p_2$ is the associative element in $\P(2)$, 
the coordinate ring of the group $\Gr_{\P}$ is generated by the 
elements of $\PP$ of the form $u\cdot\id=\gamma_{\P}(p_2;u,\id)$, where 
$u\in\P$, that is $\Hr_{\P}=\Q[\P\cdot\id]$. As a free algebra, we can 
identify it with $\Q[\P]\cong Q[\M_{\P}]/(e-1)$. 
The coproduct on $\Hr_{\P}$ which represents the composition of the series 
in $\Gr_{\P}(A)$ can be found by dualizing the inclusion of $\Gr_{\P}(A)$ 
into $\Gd_{\P}(A)$. The result makes sense on the non-commutative algebras. 

Let us then consider the surjective algebra homomorphism 
$P:\Hdnc_{\P}\longrightarrow\Hrnc_{\P}$ which sends each generator of 
$\Hdnc_{\P}=\Q\langle\PP\rangle$ of the form $u\cdot\id$ to its quotient 
$u\in\P$ and all the others to zero. 
If we verify that
\begin{align}
\label{Dr=PDd}
\Dr_{\P}(u) &= (P\otimes P)\Dd_{\P}(u\cdot\id), 
\end{align}
and that $\varepsilon(u)=P(\varepsilon_{\Hdnc_{\P}}(u\cdot\id))$, 
the coassociativity of the coproduct $\Dr_{\P}$ then follows easily. 
In fact, it holds on the commutative quotient $\Hr_{\P}$ because 
$\Gr_{\P}(A)$ is a group, and the passage to the non-commutative lift is 
as usual straightforward. All the other assertions are then easily verified. 

Therefore, it only remains to show the equality (\ref{Dr=PDd}). 
Let us fix $u\in\P$. In the sum 
\begin{align*}
\Dd_{\P}(u\cdot\id) 
&=\underset{u\cdot\id=
\gamma(\tilde{p};\tilde{q}_1,\dots,\tilde{q}_{|\tilde{p}|})}
{\underset{\tilde{q}_1,\dots,\tilde{q}_{|\tilde{p}|}\in\M_{\P}}
{\sum_{\tilde{p}\in\M_{\P}}}} \ \ 
\tilde{p} \otimes \tilde{q}_1\dots \tilde{q}_{|\tilde{p}|} ,
\end{align*}
the element $\gamma(\tilde{p};\tilde{q}_1,\dots,\tilde{q}_{|\tilde{p}|})$ can 
be of the form $u\cdot\id$ only if $\tilde{p}=p\cdot\id$
(hence $|\tilde{p}|=|p|+1$), and 
$\tilde{q}_{|\tilde{p}|}=q_{|\tilde{p}|}\cdot\id$, 
with $q_{|\tilde{p}|}\in\M_{\P}$. 
The term corresponding to $\tilde{p}=\id$, that is $p=e$, gives 
\begin{align*}
\underset{u\cdot\id=\gamma(\id;\tilde{q}_1)=\tilde{q}_1}
{\sum_{\tilde{q}_1\in \P}} \id \otimes \tilde{q}_1 
&= \id \otimes u\cdot\id. 
\end{align*}
Therefore, separating from the sum the term corresponding to
$\tilde{p}=\id$, we obtain 
\begin{align*}
\Dd_{\P}(u\cdot\id)  
&= \id \otimes u\cdot\id + \underset{u\cdot\id=
\gamma(p\cdot\id;\tilde{q}_1,\dots,\tilde{q}_{|p|},q_{|p|+1}\cdot\id)}
{\underset{\tilde{q}_1,\dots,\tilde{q}_{|p|}\in \P,\ q_{|p|+1}\in \M_{\P}}
{\sum_{p\in \P}}} \ \ 
p\cdot\id \otimes \tilde{q}_1\dots \tilde{q}_{|p|}\ (q_{|p|+1}\cdot\id). 
\end{align*}
Since for any $p\in \P$ we have $p\cdot\id = \gamma(p_2;p,\id)$, using 
the associativity of $\gamma$ we obtain 
\begin{align*}
u\cdot\id &= \gamma(\gamma(p_2;p,\id);
\tilde{q}_1,\dots,\tilde{q}_{|p|},q_{|p|+1}\cdot\id) \\ 
&= \gamma(p_2;\gamma(p;\tilde{q}_1,\dots,\tilde{q}_{|p|}),
\gamma(\id;q_{|p|+1}\cdot\id)) \\ 
&= \gamma(p;\tilde{q}_1,\dots,\tilde{q}_{|p|}) \cdot q_{|p|+1}\cdot \id. 
\end{align*}
Therefore $u=\gamma(p;\tilde{q}_1,\dots,\tilde{q}_{|p|})\cdot q_{|p|+1}$, 
and we get 
\begin{align*}
\Dd_{\P}(u\cdot\id)  
&= \id \otimes u\cdot\id +
\underset{u=\gamma(p;\tilde{q}_1,\dots,\tilde{q}_{|p|})\cdot q_{|p|+1}}
{\underset{\tilde{q}_1,\dots,\tilde{q}_{|p|}\in \P,\ q_{|p|+1}\in \P}
{\sum_{p\in \P}}} \ \  
p\cdot \id \otimes \tilde{q}_1\dots \tilde{q}_{|p|}\ (q_{|p|+1}\cdot\id). 
\end{align*}
Now, applying the projection $P$, we kill all the elements $\tilde{q}_i$ of
$\Hdnc_{\P}$ different from $\tilde{q}_i=q_i\cdot\id$, therefore the sum on 
the right hand-side is reduced to 
\begin{align*}
(P\otimes P)\Dd_{\P}(u\cdot\id)  
&= (P\otimes P)\big[\id \otimes u\cdot\id + 
\underset{u=\gamma(p;q_1\cdot\id,\dots,q_{|p|}\cdot\id)\cdot q_{|p|+1}}
{\underset{q_1,\dots,q_{|p|},q_{|p|+1}\in \M_{\P}} {\sum_{p\in \P}}} \ \  
p\cdot \id \otimes (q_1\cdot\id)\cdots (q_{|p|+1}\cdot\id)\big]
\end{align*} 
and we finally obtain 
\begin{align*}
(P\otimes P)\Dd_{\P}(u\cdot\id)  
&= 1 \otimes u 
+ \underset{u=\gamma(p;q_1\cdot\id,\dots,q_{|p|}\cdot\id)\cdot q_{|p|+1}}
{\underset{q_1,\dots,q_{|p|+1}\in \M_{\P}} {\sum_{p\in \P}}} \ \  
p \otimes q_1\cdots q_{|p|+1}. 
\end{align*}
\end{proof}

Note that the coordinate ring of the group $\Gr_{\P}$ has the same generators 
as the coordinate ring of the group $\Gi_{\M_{\P}}$, and the same holds for  
their non-commutative lifts. As algebras they are isomorphic, but they 
differ for the coalgebra structure. 
Since the action of $\Gd_{\P}$ on $\Gr_{\P}$ is induced by that on 
$\Gi_{\M_{\P}}$, it is not surprising that the dual coactions coincide. 

In other words, if we define the algebra morphism 
$\ddif_{\P}: \Hrnc_{\P}\longrightarrow \Hrnc_{\P}\otimes\Hdnc_{\P}$ 
by setting, on the generators $u\in\P$, 
\begin{align}
\label{ddif_P-def}
\ddif_{\P}(u) &:= \underset{\gamma(p;q_1,\dots,q_{|p|})=u}
{\underset{q_1,...,q_{|p|}\in\PP\cup\{\id\}}{\sum_{p\in\P}}} 
p \otimes q_1\dots q_{|p|}, 
\end{align}
then the following result is straightforward. 

\begin{proposition}
\label{Hr-coaction}
The map $\ddif_{\P}$ is a graded right coaction of the Hopf algebra 
$\Hdnc_{\P}$ on the algebra $\Hrnc_{\P}$. 
\end{proposition}


\subsection{Diffeomorphisms subgroups of tree-expanded series}

Let us apply the results of the previous section to the duplicial operad.
As observed in Section~\ref{subsection-action-tree}, the operad $\Dup$ 
has two associative binary operations, $\under$ and $\over$, which lead 
to the two groups $\Gu(A)$ and $\Go(A)$ of invertible tree-expanded series. 
Each of the two operations determines two linear maps 
$\lambda,\rho:A[[Y]]\longrightarrow A[[\Y]]$ and consequently two subgroups 
$\Gl_Y(A)$ and $\Gr_Y(A)$ of $\Gd_Y(A)$. 

For our purposes, we are only interested in one of the four resulting groups: 
the group $\Gr_Y(A)$ corresponding to the operation $\over$. However, 
in order to discuss some of its properties, we also make use of the group 
$\Gl_Y(A)$ corresponding to the operation $\under$. 
To fix the notations, we recall these two groups explicitly: 
\begin{align*}
\Gr_Y(A) &:= \Go(A)\over x^{\treeA} 
= \big\{ \rho_f(x)= \sum_{t\in Y} f_t \ x^{\lvertexgraft{t}}, 
\quad f_t  \in  A,\ f_{\treeO}=1 \big\}, \\  
\Gl_Y(A) &:= x^{\treeA}\under\Gu(A) 
= \big\{ \lambda_f(x)= \sum_{t\in Y} f_t \ x^{\rvertexgraft{t}}, 
\quad f_t \in A,\ f_{\treeO}=1 \big\}. 
\end{align*}
The intersection $\Gr_Y(A)\cap\Gl_Y(A)$ obviously contains only the unit\ 
$\id(x)$.

As we already observed, in the case $\P=\As$ all these subgroups in fact 
coincide with the whole group of formal diffeomorphisms. In the present case 
this surely does not hold. Moreover the two subgroups $\Gr_Y(A)$ and 
$\Gl_Y(A)$ are not normal in $\Gd_Y(A)$, however they allow to reconstruct 
the group $\Gd_Y(A)$. To do this, let us fix the notation 
\begin{align*}
\Gl_Y(A)\circ\Gr_Y(A) & := \big\{ \lambda_f \circ \rho_g \ 
\mbox{where $\lambda_f \in\Gl_Y(A)$ and $\rho_g \in\Gr_Y(A)$} \big\}. 
\end{align*}

\begin{lemma}
Each series in $\Gd_Y(A)$ can be written as the composition of two series 
in $\Gl_Y(A)$ and $\Gr_Y(A)$, that is 
\begin{align}
\label{composition-decomposition}
\Gd_Y(A) &= \Gl_Y(A)\circ\Gr_Y(A) = \Gr_Y(A)\circ\Gl_Y(A). 
\end{align} 
Moreover, this decomposition is unique if, on the left hand-side, we restrict 
the choice to the comb-trees, that is 
\begin{align}
\label{composition-decomposition-unique}
\Gd_Y(A) & = \iu(\Gd(A))\circ \Gr_Y(A)= \io(\Gd(A)) \circ\Gl_Y(A).  
\end{align} 
\end{lemma}

\begin{proof} 
To show the equality (\ref{composition-decomposition}), we have
to show that any tree-expanded formal diffeomorphism
$\eta(x)=\sum_{u\in \Y} \eta(u)\ x^u$ can be written as the
compositions $\lambda_f \circ\rho_g$ and $\rho_{g'}\circ\lambda_{f'}$, 
for some $f,g,f',g' \in \Gi_Y(A)$.  In other words, since the coefficients 
lie in a commutative unital algebra and can be chosen arbitrarily, 
we have to show that, in the compositions
$\lambda_f \circ\rho_g$ and $\rho_{g'}\circ\lambda_{f'}$, the power
$x^u$ appears for all the trees $u\in \Y$.  Let us show it for the
case $\lambda_f \circ\rho_g$, the same procedure can be adapted to
the other case.

In the composition of the two series 
$\lambda_f(x)=\sum_{t\in Y} f_t x^{\treeAsmall\under t}$ and 
$\rho_g(x)=\sum_{s\in Y} g_s x^{s\over\treeAsmall}$, namely
\begin{align*}
\big(\lambda_f\circ\rho_g\big)(x) 
&=\sum_{t,s_0,s_1,\dots,s_{|t|}\in Y} f_t\  g_{s_0} g_{s_1}\dots g_{s_{|t|}}\ 
x^{\mu_{\treeAsmall \under t}
(s_0\over\treeAsmall,\dots,s_{|t|}\over\treeAsmall)}, 
\end{align*}
there appears the power $x^u$ for 
$u=\mu_{\treeA\under t}(s_0 \over\treeA,\dots,s_{|t|} \over \treeA)$, 
where $t$ and $s_0,\dots,s_{|t|}$ are arbitrary trees (all possibly equal to
$\treeO$).

If $t=\treeO$, we get
\begin{equation*}
u=\mu_{\treeAsmall}(s_0 \over \treeA) = s_0 \over \treeA = \lvertexgraft{s_0}. 
\end{equation*}
Since $s_0$ runs over all possible trees, this $u$ recovers all
trees with nothing branched on the right of the root.

If $t\neq \treeO$, we use the fact that 
$\treeA\under t =\mu_{\treeBA}(\treeA,t)$, the associativity of the product 
$\mu$ shaped by trees, and the associativity of the over and under
products, to get
\begin{align}
u &= \mu_{\mu_{\treeBA}(\treeA,t)}
(s_0 \over\treeA,\dots,s_{|t|}\over\treeA) \nonumber \\ 
&= \mu_{\treeBA}\left(\mu_{\treeA}(s_0\over\treeA),
\mu_t(s_1\over\treeA,\dots,s_{|t|}\over\treeA)\right) \nonumber \\ 
&= s_0\over\treeA\under\mu_t(s_1\over\treeA,\dots,s_{|t|}\over\treeA). 
\label{decomposition-formula}
\end{align}
Since $t$ and $s_0,\dots s_{|t|}$ run over all possible trees
(including the root-tree $\treeO$ for the $s_i$'s), we can recover
any possible tree $u\in \Y$ with something branched simultaneously on
the left and on the right of the root.
  
The above decomposition is clearly not unique, because different choices 
of $t$ and $s_1,\dots s_{|t|}$ might give rise to the same tree $u$. 
For instance, if in Eq.~(\ref{decomposition-formula}) we choose 
$t=\treeAB$, any $s_1$, and $s_2=\treeO$, we get
\begin{align*}
\mu_t(s_1 \over \treeA,s_2 \over \treeA) &= 
\mu_{\treeAB}(\lvertexgraft{s_1},\treeA) = 
\setlength{\unitlength}{4pt}
\psset{unit=4pt}
\psset{runit=2.8pt}
\psset{linewidth=0.2}
\begin{pspicture}(0,0)(8,6)
\psline(5,-.5)(5,1)
\psline(5,1)(3,3)
\psline(5,1)(6,2)
\psline(4,2)(5,3)
\put(1,5){$s_1$}
\end{pspicture}
\end{align*}
and therefore for any $s_0$ we get 
$u= s_0 \over\ \treeA \under
  \setlength{\unitlength}{4pt} \psset{unit=4pt} \psset{runit=2.8pt}
  \psset{linewidth=0.2}
  \begin{pspicture}(0,0)(8,6)
    \psline(5,-.5)(5,1)
    \psline(5,1)(3,3)
    \psline(5,1)(6,2)
    \psline(4,2)(5,3)
    \put(1,5){$s_1$}
  \end{pspicture} = 
  \setlength{\unitlength}{5pt}
  \psset{unit=4pt}
  \psset{runit=2.8pt}
  \psset{linewidth=0.2}
  \begin{pspicture}(0,0)(7,8)
    \psline(5,-.5)(5,1)
    \psline(5,1)(7,3)
    \psline(5,1)(4,2)
    \psline(6,2)(4,4)
    \psline(5,3)(6,4)
    \put(1.5,3){$s_0$}
    \put(2,7){$s_1$}
  \end{pspicture}$. 
But if we choose $t'=\treeA$ and $s'_1=\lvertexgraft{s_1}$, we
get the same result for any $s_0$ because
\begin{align*}
    \mu_{t'}(s'_1 \over \treeA) &=
    \mu_{\treeA}(\setlength{\unitlength}{4pt} \psset{unit=4pt}
    \psset{runit=2.8pt} \psset{linewidth=0.2}
    \begin{pspicture}(0,0)(8,6)
      \psline(5,-.5)(5,1)
      \psline(5,1)(3,3)
      \psline(5,1)(6,2)
      \psline(4,2)(5,3)
      \put(1,5){$s_1$}
    \end{pspicture}) 
    = 
    \setlength{\unitlength}{4pt}
    \psset{unit=4pt}
    \psset{runit=2.8pt}
    \psset{linewidth=0.2}
    \begin{pspicture}(0,0)(8,6)
      \psline(5,-.5)(5,1)
      \psline(5,1)(3,3)
      \psline(5,1)(6,2)
      \psline(4,2)(5,3)
      \put(1,5){$s_1$}
    \end{pspicture}. 
\end{align*}

To show that the decomposition (\ref{composition-decomposition-unique}) 
is unique, in Eq.~(\ref{decomposition-formula}) it suffices to consider, 
for $t$, only the right-comb trees $\combR$. With this choice, we get
\begin{equation*}
u=\combRgraft{s_0}{s_1}{s_{|t|}}, 
\end{equation*}
and therefore, for arbitrary $s_0,\dots,s_{|t|} \in Y$, with
$|t|\geq 1$, we recover in a unique way all trees $u\in \Y$ with
something branched simultaneously on the left and on the right of the root.
Then we apply Proposition~\ref{Gd_P->Gd} to identify the group $\Gd(A)$ 
with the subgroup of $\Gd_Y(A)$ made of series expanded only on the 
right-comb trees. 
We have therefore proved the uniqueness of decomposition 
$\iu(\Gd(A))\circ\Gr_Y(A)=\Gd_Y(A)$.
The same argument applies to $\io(\Gd(A))\circ\Gl_Y(A)=\Gd_Y(A)$.
\end{proof}

The order map $\pi$ also gives two surjective group morphisms from
$\Gr_Y(A)$ and $\Gl_Y(A)$ to $\Gd(A)$.
In fact, since $\Gr_Y(A)$ and $\Gl_Y(A)$ are subgroups of $\Gd_Y(A)$, 
and $\pi$ is a group homomorphism from $\Gd_Y(A)$ to $\Gd(A)$, 
it only remains to show that $\pi$ is still surjective when restricted to 
$\Gr_Y(A)$ or $\Gl_Y(A)$. This follows from the fact that $\Gr_Y(A)$ and 
$\Gl_Y(A)$ contain $\Gd(A)$ via the inclusions of Proposition~\ref{Gd_P->Gd}, 
which are sections of $\pi$.
\bigskip 

To conclude this section, we apply Theorem~\ref{Hr_P} to describe 
explicitly the Hopf structure of the algebra $\Hr_Y$, because it gives rise 
to another new Hopf algebra on trees, which is neither commutative nor 
cocommutative. 

\begin{corollary}
The free associative algebra 
$\Hrnc_Y=\Q\langle\Y\rangle\cong\Q\langle Y\rangle/(\treeO-1)$ 
is a graded and connected algebra, with grading given by the 
order of trees. The coproduct $\Dr_Y$ is defined on any $u\in\Y$ by 
\begin{align}
\Dr_Y(u) &:= 1 \otimes u 
+ \underset{u=\mu_t(s_1\over\treeAsmall,\dots,s_{|t|}\over\treeAsmall)
\over s_{|t|+1}}
{\underset{s_1,\dots,s_{|t|+1}\in Y} {\sum_{t\in \Y}}} \ \  
t \otimes s_1\dots s_{|t|} s_{|t|+1}. 
\label{Dr_Y(u)}
\end{align}
and the counit is $\varepsilon(u)=0$ for any $u\in \Y$.
\end{corollary}

For instance, setting $\treeO=1$, the coproduct on small trees is: 
\begin{align*}
\Dr_Y(\treeA) &= \treeA \otimes 1 + 1 \otimes \treeA, \\ 
\Dr_Y(\treeAB) &= \treeAB \otimes 1 + 2 \treeA \otimes \treeA 
+ 1 \otimes \treeAB, \\
\Dr_Y(\treeBA) &= \treeBA \otimes 1 + 1 \otimes \treeBA, \\
\Dr_Y(\treeABC) &= \treeABC \otimes 1 + 3 \treeAB \otimes \treeA 
+ \treeA \otimes (2\treeAB +\treeA^2) + 1 \otimes \treeABC, \\
\Dr_Y(\treeBAC) &= \treeBAC \otimes 1 + \treeBA \otimes \treeA
+ \treeA \otimes \treeBA + 1 \otimes \treeBAC, \\
\Dr_Y(\treeACA) &= \treeACA \otimes 1 + \treeBA \otimes \treeA 
+ \treeA \otimes \treeBA + 1 \otimes \treeACA, \\
\Dr_Y(\treeCAB) &= \treeCAB \otimes 1 + \treeBA \otimes \treeA 
+ 1 \otimes \treeCAB, \\
\Dr_Y(\treeCBA) &= \treeCBA \otimes 1 + 1 \otimes \treeCBA. 
\end{align*} 


\subsection{Subgroup dual to the Hopf algebra $\Ha$}
\label{subsection-charge}

The main aim of this section is to define the subgroup $\Ga$  
represented by the ``charge Hopf algebra'' $\Ha$ introduced in
\cite{BFqedtree}, and used in \cite{BFqedren} to describe the 
renormalization of the electric charge in quantum electrodynamics. 
\bigskip 

Let us fix an associative, commutative and unital algebra $A$. 
For any $f(x)=\sum_{t\in Y} f_t\ x^t\in A[[Y]]$, the series 
$x^{\treeO}-x^{\treeA}\under f(x)$ belongs to the set $\Gi_Y(A)$, 
and therefore to the group $\Go(A)$ of tree-expanded invertible series 
with respect to the product $\over$. 
Let us call $\tilde{f}(x)=(x^{\treeO}-x^{\treeA}\under f(x))^{-1}$ 
its inverse in $\Go(A)$, and set 
\begin{align}
\label{Ga}
\Ga(A) &:= \left\{ \alpha_f(x)=\rho_{\tilde{f}}(x)
=(x^{\treeO}-x^{\treeA}\under f(x))^{-1}\over x^{\treeA},\ 
f\in A[[Y]] \right\}. 
\end{align}
A tree-expanded diffeomorphism $\alpha_f(x)$ can be thought as the 
translations by the series $f$ which fixes zero. 

For any tree $t\in Y$, set $V(t)=\treeA\under t$. 

\begin{lemma}
\label{Ga-Gr}
The set $\Ga(A)$ coincides with the subset of $\Gr_Y(A)$ made of the series 
$\rho_g(x)=\sum_{t\in Y} g_t\ x^{\lvertexgraft{t}}$ such that 
\begin{align}
\label{Gr-Ga}
g_t &= g_{V(t_1)}\ g_{V(t_2)}\ \cdots\ g_{V(t_n)}, 
\qquad\mbox{if\quad $t=V(t_1)\over V(t_2)\over\cdots\over V(t_n)$}. 
\end{align}
\end{lemma} 

\begin{proof}
For any fixed $f(x)\in A[[Y]]$, the inverse of the series 
$x^{\treeO}-x^{\treeA}\under f(x)= 
x^{\treeO}-\sum_{t\in Y} f_t\ x^{\rvertexgraft{t}}$ in $\Go(A)$ 
is 
\begin{align*}
\tilde{f}(x) &= x^{\treeO} 
+ \sum_{n=1}^\infty \left(\sum_{t\in Y} 
f_t\ x^{\rvertexgraft{t}}\right)^{\over n} 
= x^{\treeO} 
+ \sum_{n=1}^\infty \sum_{t_1,...,t_n\in Y} f_{t_1}\cdots f_{t_n}\ 
x^{V(t_1)\over\cdots\over V(t_n)}. 
\end{align*}
Any tree $t \neq \treeO$ can be written in a unique way as an over product of 
trees which have nothing branched at the left of the root, in fact
\begin{align*}
t &= \combLgraft{t_1}{t_2}{t_n}=V(t_1)\over V(t_2)\over\cdots\over V(t_n).  
\end{align*}
Therefore we have 
\begin{align*}
\tilde{f}(x) &= x^{\treeO} 
+ \sum_{n=1}^\infty \underset{t=V(t_1)\over\cdots\over V(t_n)}{\sum_{t\in Y}} 
f_{t_1}\cdots f_{t_n}\ x^t . 
\end{align*}

Then, varying $f\in A[[Y]]$, the series 
$\alpha_f(x)=\tilde{f}(x)\over x^{\treeA}$ give exactely all the series 
$\rho_g(x)$, where $g(x)=x^{\treeO}+\sum_{t\in\Y} g_t\ x^t$ has 
arbitrary coefficients $g_{V(t)}= f_t$ and constrained coefficients 
$g_{V(t_1)\over V(t_2)\over\cdots\over V(t_n)} 
= g_{V(t_1)}\ g_{V(t_2)}\ \cdots\ g_{V(t_n)}$. 
\end{proof}

\begin{theorem}
The set $\Ga(A)$ is a subgroup of $\Gr_Y(A)$. 
\end{theorem} 

\begin{proof}
Since the series $f(x)=0$ gives $\alpha_f(x)=x^{\treeA}$, it suffices 
to show that the subset $\Ga(A)$ is closed for the composition. 

Let us exploit Lemma~\ref{Ga-Gr}, and choose two generic series in $\Ga(A)$ 
by taking two series $\rho_f$ and $\rho_g$ in $\Gr_Y(A)$ such that the 
coefficients of the series $f,g\in\Go(A)$ satisfy the condition~(\ref{Gr-Ga}). 
Let $h\in\Go(A)$ be the series which results from the composition 
$\rho_f\circ\rho_g=\rho_h$. We have to show that the coefficients of $h$ 
also satisfy the condition~(\ref{Gr-Ga}), that is, we have to show that 
for any $u\in\Y$ we have $h_u=h_{V(u_1)}\cdots h_{V(u_n)}$, if 
$u=V(u_1)\over\cdots\over V(u_n)$. For this, it suffices to show that 
\begin{align}
\label{h_ul/V(ur)}
h_u=h_{u^l} h_{V(u^r)} \quad\mbox{if $u=u^l\over V(u_r)$}. 
\end{align}

Applying the definition of the composition, for any $u\in\Y$ we have 
\begin{align}
\label{h_u}
h_u &= g_u + 
\sum_{u=\mu_t(s_1\over\treeAsmall,...,s_{|t|}\over\treeAsmall)\over s_{|t|+1}}
f_t\ g_{s_1}\cdots g_{s_{|t|+1}}, 
\end{align}
where from now on we suppose that the sums run over all the trees in the 
set $\Y$ if they appear as subindeces of the tree-product $\mu$ 
(in this case $t$), and to the set $Y$ if they appear inside the arguments 
of $\mu$ or anywhere else (in this case $s_1,...,s_{|t|+1}$). 

In particular, we need an explicit expression the coefficient $h_{V(u)}$, 
where $u\in Y$. Let us compute it. 
If $u=\treeO$ and $V(u)=\treeA$, it is easy to see that 
\begin{align*}
h_{\treeAsmall} &= g_{\treeAsmall} + f_{\treeAsmall}. 
\end{align*}
Then we suppose that $u\neq\treeO$. 
In Eq.~(\ref{h_u}), we replace the tree $u$ by the tree 
$V(u)=\rvertexgraft{u}$, and obtain
\begin{align*}
h_{V(u)} &= g_{V(u)} + \sum_{V(u)=\mu_{\tilde{t}}
(\tilde{s}_1\over\treeAsmall,\dots,\tilde{s}_{|\tilde{t}|}\over\treeAsmall)
\over \tilde{s}_{|\tilde{t}|+1}}
f_{\tilde{t}}\ g_{\tilde{s}_1}\cdots g_{\tilde{s}_{|\tilde{t}|+1}}.  
\end{align*}
The tree
$\mu_{\tilde{t}}(\tilde{s}_1\over\treeA,\dots,\tilde{s}_{|\tilde{t}|}
\over\treeA)\over \tilde{s}_{|\tilde{t}|+1}$ can be of the form
$V(u)=\rvertexgraft{u}$ only if $\tilde{s}_{|\tilde{t}|+1}=\treeO$,
$\tilde{t}= V(t)=\rvertexgraft{t}$ with $t\in Y$, and $\tilde{s}_1=\treeO$.  
The case $t=\treeO$ corresponds to $u=\treeO$, that we already computed apart. 
For $t\neq\treeO$, we write $V(t)=\mu_{\treeBA}(\treeA,t)$ and apply the 
associativity of $\mu$ to conclude that
\begin{align*}
V(u)&=\mu_{V(t)}
(\treeA,\tilde{s}_2\over\treeA,\dots,\tilde{s}_{|t|+1}\over\treeA)
\over\treeO \\
&= V(\mu_t(\tilde{s}_2\over\treeA,\dots,\tilde{s}_{|t|+1}\over\treeA)), 
\end{align*}
and therefore $u=\mu_t(\tilde{s}_2\over\treeA,\dots,\tilde{s}_{|t|+1}\treeA)$. 
By renaming the trees $\tilde{s}_i=s_{i-1}$, we finally obtain 
\begin{align}
h_{V(u)}
&= g_{V(u)} 
+ \sum_{u=\mu_t(s_1\over\treeAsmall,\dots,s_{|t|}\over\treeAsmall)}
f_{V(t)}\ g_{s_1}\cdots g_{s_{|t|}}. 
\label{h_{V(u)}}
\end{align}

Let us now prove (\ref{h_ul/V(ur)}). We start again from (\ref{h_u}), 
for a fixed tree $u=u^l\over V(u^r)$.
\bigskip  

\noindent{\em Assume $u^r=\treeO$.\/} 
Let us start by considering the case $u=u^l\over\treeA$. 
We already computed $h_{\treeA}=g_{\treeA}+f_{\treeA}$. 
The sum in (\ref{h_u}) is over all trees $t\in \Y$ and 
$s_1,\dots,s_{|t|+1}\in Y$ such that 
$u=u^l\over\treeA=\mu_t(s_1\over\treeA,\dots,s_{|t|}\over\treeA)\over
s_{|t|+1}$.  
Let us list the contributions to this sum coming from different cases. 
\bigskip 

\noindent{\em Case 1.\/} 
If $s_{|t|+1}=\treeO$, the equality 
$u=u^l\over\treeA=\mu_t(s_1\over\treeA,\dots,s_{|t|}\over\treeA)$ 
is possible if and only if $t=t^l\over\treeA$ and 
$s_{|t|}=s_{|t|}^l\over\treeA$. 
Then we distinguish the following two possible cases. 
\bigskip 

\noindent{\em Case 1a.\/} 
If $t^l\neq\treeO$, and therefore $|t|=|t^l|+1$, then 
$u^l=\mu_{t^l}(s_1\over\treeA,\dots,s_{|t^l|}\over\treeA)\over s_{|t^l|+1}^l$. 
We then rename $t^l=:w$. Since $g_{\treeO}=1$, and 
$f_{w\over\treeAsmall}= f_{w} f_{\treeAsmall}$, we have the contribution 
\begin{align*}
\sum_{u^l=\mu_{w}(s_1\over\treeAsmall,\dots,s_{|w|}\over\treeAsmall)
\over s_{|w|+1}}
f_{w} f_{\treeAsmall}\ g_{s_1}\cdots g_{s_{|w|+1}}. 
\end{align*}
\bigskip 

\noindent{\em Case 1b.\/} 
If $t^l=\treeO$, then 
$u=u^l\over\treeA=\mu_{\treeAsmall}(s_1\over\treeA)=s_1\over\treeA$, 
and therefore $s_1=u^l$. We then have the contribution 
\begin{align*}
&f_{\treeAsmall}\ g_{u^l}. 
\end{align*}
\bigskip 

\noindent{\em Case 2.\/} 
If $s_{|t|+1}\neq\treeO$, the equality 
\begin{align*}
u^l\over\treeA &= \mu_t(s_1\over\treeA,\dots,s_{|t|}\over\treeA)
\over s_{|t|+1}
\end{align*}
is possible if and only if $s_{|t|+1}=s_{|t|+1}^l\over\treeA$ and 
$\mu_t(s_1\over\treeA,\dots,s_{|t|}\over\treeA)\over s_{|t|+1}^l=u ^l$.  
We then rename the free trees $t=:w$ and $s_{|t|+1}^l=:s_{|t|+1}$, 
and obtain the contribution
\begin{align*}
&\sum_{u^l=\mu_w(s_1\over\treeAsmall,\dots,s_{|w|}\over\treeAsmall)
\over s_{|w|+1}}
f_w\  g_{s_1}\cdots g_{s_{|w|+1}} g_{\treeAsmall}. 
\end{align*}
\bigskip 

Summing up all the possible cases, we obtain 
\begin{align*}
h_{u^l\over\treeAsmall}   
&= \left( g_{u^l}
+\sum_{u^l=\mu_w(s_1\over\treeAsmall,\dots,s_{|w|}\over\treeAsmall)
\over s_{|w|+1}}
f_w\  g_{s_1}\cdots g_{s_{|w|+1}}\right)
\left( g_{\treeAsmall}+f_{\treeAsmall}\right) 
= h_{u^l} h_{\treeAsmall}. 
\end{align*}
\bigskip 

\noindent{\em Assume $u^r \neq \treeO$.\/} 
The sum in (\ref{h_u}) is over all trees $t\in \Y$ and 
$s_1,\dots,s_{|t|+1}\in Y$ such that 
$u=u^l\over V(u^r)=\mu_t(s_1\over\treeA,\dots,s_{|t|}\over\treeA)\over
s_{|t|+1}$.  
Let us list all the contributions to this sum. 
\bigskip 

\noindent{\em Case 1.\/} 
If $s_{|t|+1}=\treeO$, since $t \neq\treeO$, we split $t=t^l\over V(t^r)$,
and distinguish the following three possible cases. 
\bigskip 

\noindent{\em Case 1a.\/} 
If $t^l\neq\treeO$ and $t^r =\treeO$, then $|t|=|t^l|+1$ and 
the equality becomes 
\begin{align*}
u^l\over V(u^r) 
&= \mu_{t^l\over\treeA}(s_1\over\treeA,\dots,s_{|t^l|+1}\over\treeA) 
= \mu_{t^l}(s_1\over\treeA,\dots,s_{|t^l|}\over\treeA) 
\over s_{|t^l|+1}\over\treeA \\
& = \mu_{t^l}(s_1\over\treeA,\dots,s_{|t^l|}\over\treeA) 
\over s_{|t^l|+1}\over V(\treeO), 
\end{align*}
which is impossible because $u^r\neq\treeO$. 
\bigskip 

\noindent{\em Case 1b.\/} 
If $t^l=\treeO$ and $t^r\neq\treeO$, then $t=V(t^r)=\mu_{\treeBA}(\treeA,t^r)$ 
with $|t|=|t^r|+1$. 
Using the associativity of the $\mu$ product we get
\begin{align*}
\mu_{V(t^r)}(s_1\over\treeA,\dots,s_{|t^r|+1}\over\treeA) &= 
s_1\over\treeA\under \mu_{t^r}(s_2\over\treeA,\dots,s_{|t^r|+1}\over\treeA). 
\end{align*}
This tree can be equal to $u=u^l\over\treeA\under u^r$ if and only if $s_1=u^l$
and $\mu_{t^r}(s_2\over\treeA,\dots,s_{|t^r|+1}\over\treeA)=u^r$. 
Let us rename the free trees as follows:  
$t^r=:y$, $s_i=:z_{i-1}$ for $i=2,\dots,|t^r|+1$. Then the contribution 
to the whole sum for this case is 
\begin{align*}
& \sum_{u^r=\mu_y(z_1\over\treeAsmall,\dots,z_{|y|}\over\treeAsmall)}
\!\!\!\!\!\!\!\!\!\!\!\! 
f_{V(y)}\ g_{u^l} g_{z_1}\cdots g_{z_{|y|}}. 
\end{align*}
\bigskip 

\noindent{\em Case 1c.\/} 
If $t^l\neq\treeO$ and $t^r\neq\treeO$, using again the properties of the
$\mu$ product we have  
\begin{align*}
u^l\over V(u^r) 
&= \mu_{t^l}(s_1\over\treeA,\dots,s_{|t^l|}\over\treeA) \over 
\mu_{V(t^r)}(s_{|t^l|+1}\over\treeA,\dots,s_{|t^l|+|t^r|+1}\over\treeA) \\ 
&= \mu_{t^l}(s_1\over\treeA,\dots,s_{|t^l|}\over\treeA) \over s_{|t^l|+1}\over
V\big(\mu_{t^r}(s_{|t^l|+2}\over\treeA,\dots,s_{|t^l|+|t^r|+1}\over\treeA)\big).
\end{align*}
We rename the free trees as $t^l=:w$, $s_i=:v_i$ for $i=1,\dots,|t^l|+1$,
and $t^r=:y$, $s_{|t^l|+j}=:z_{j-1}$ for $j=2,\dots,|t^r|+1$. 
Since $f_{w\over V(y)}=f_w\ f_{V(y)}$, we obtain the contribution  
\begin{align*}
\underset{u^r=\mu_y(z_1\over\treeAsmall,\dots,z_{|y|}\over\treeAsmall)}
{\sum_{u^l=\mu_w(v_1\over\treeAsmall,\dots,v_{|w|}\over\treeAsmall)
\over v_{|w|+1}}}
f_w f_{V(y)}\ g_{v_1}\cdots g_{v_{|w|+1}}\ g_{z_1}\cdots g_{z_{|y|}}. 
\end{align*}
\bigskip 

\noindent{\em Case 2.\/} 
If $s_{|t|+1}\neq\treeO$, then it can be decomposed as 
$s_{|t|+1}=s_{|t|+1}^l\over V(s_{|t|+1}^r)$, and therefore the
equality 
\begin{align*}
u^l\over V(u^r) &= \mu_t(s_1\over\treeA,\dots,s_{|t|}\over\treeA)
\over s_{|t|+1}^l \over V(s_{|t|+1}^r)
\end{align*}
is possible if and only if 
$\mu_t(s_1\over\treeA,\dots,s_{|t|}\over\treeA)\over s_{|t|+1}^l=u ^l$ 
and $s_{|t|+1}^r=u^r$. 
We then rename the free trees as $t=:w$, $s_i=:v_i$ for $i=1,\dots,|t|$ and 
$s_{|t|+1}^l=:v_{|w|+1}$, and obtain the contribution
\begin{align*}
&\sum_{u^l=\mu_w(v_1\over\treeAsmall,\dots,v_{|w|}\over\treeAsmall)
\over v_{|w|+1}}
f_w\ g_{v_1}\cdots g_{v_{|w|+1}} g_{V(u^r)}. 
\end{align*}
\bigskip 

Summing up all the possible cases, and rearranging them using 
(\ref{h_u}) and (\ref{h_{V(u)}}), we finally obtain 
\begin{align*}
h_{u^l\over V(u^r)} 
&= g_u 
+ \sum_{u=\mu_t(s_1\over\treeAsmall,\dots,s_{|t|}\over\treeAsmall)
\over s_{|t|+1}}
f_t\ g_{s_1}\cdots g_{s_{|t|+1}} \\ 
& = \left( g_{u^l} + 
\sum_{u^l=\mu_w(v_1\over\treeAsmall,\dots,v_{|w|}\over\treeAsmall)
\over v_{|w|+1}}
f_w\ g_{v_1}\cdots g_{v_{|w|+1}} \right) \times \\ 
& \hspace{5cm} \left( g_{V(u^r)} + 
\sum_{u^r=\mu_y(z_1\over\treeAsmall,\dots,z_{|y|}\over\treeAsmall)}
f_{V(y)}\ g_{z_1}\cdots g_{z_{|y|}} \right) \\
&= h_{u^l}\ h_{V(u^r)}. 
\end{align*}
\end{proof}

The construction of the group $\Ga(A)$ is clearly functorial in $A$. 

\begin{theorem}
The group functor $\Ga$ is represented by the Hopf algebra $\Ha$. 
\end{theorem} 

\begin{proof}
The functoriality of the group $\Ga(A)$, with respect to the algebra $A$, 
is obvious, as well as the fact that $\Ga$ is a proalgebraic group. 
In fact, after Lemma~\ref{Ga-Gr}, the coordinate ring of $\Ga$ is the 
quotient of $\Hr_Y$ by the ideal generated by the relation 
$t=V(t_1)\over V(t_2)\over\cdots\over V(t_n)$, for all $t\in\Y$, 
and therefore it is freely spanned by the trees of the form 
$V(t)\over\treeA$, for any $t\in Y$. 
In other words, the coordinate ring of $\Ga$ is the polynomial algebra 
$\Q[V(t)\over\treeA,t\in Y]$, which is isomorphic, as an algebra, to the 
polynomial algebra $\Q[V(t),t\in Y]$. 

Let us recall, from \cite{BFqedtree}, that $\Ha$ is the abelian quotient 
of the algebra $\Q Y$ of all trees endowed with the over product. 
Thus the root tree $\treeO$ is the unit, and the algebra $\Ha$ is in fact 
isomorphic to the polynomial algebra $\Q[V(t), t\in Y]$. 
In \cite{BFqedtree} it was shown that $\Ha$ is a connected graded Hopf
algebra, with the grading given by the order of the trees. 
The coproduct $\Da:\Ha \longrightarrow \Ha \otimes\Ha$ is the algebra
morphism defined on the generators by the assignment 
\begin{align}
\label{Da-def}
\Da (V(t)) &= 1\otimes V(t)+\da(V(t)), 
\end{align}
where $\da:\Ha\longrightarrow \Ha\otimes\Ha$ is a right coaction of $\Ha$ 
on itself (w.r.t. the coproduct $\Da$), defined recursively as 
\begin{align}
\label{da-def}
\da(V(t)) &= (V\otimes\Id) \left[\Da(t^l)\over\da(V(t^r))\right],
\end{align} 
where $t=t^l\over V(t^r)$. 
The counit $\varepsilon:\Ha\longrightarrow \Q$ is the algebra 
morphism with value $\varepsilon(V(t))=0$ on the generators. 

It is clear that the coordinate ring of $\Ga$ is isomorphic to $\Ha$, 
as an algebra. It remains to show that the coproduct dual to the composition 
is indeed $\Da$. To do this, we prove that the projection 
$R:\Hr_Y=\Q[\Y] \longrightarrow \Ha\cong\Q[V(t),t\in Y]$,  
dual to the inclusion of $\Ga(A)$ into $\Gr_Y(A)$, is a morphism of Hopf 
algebras, that is 
\begin{align}
\label{Dr-Da}
\Da(R(u))=(R\otimes R)\Dr_Y(u)
\end{align}
for all $u\in\Y$.
The map $R$ is the algebra morphism which sends the generators $u\in\Y$ of 
$\Hr_Y$ into themselves, seen as over products of its components, that is 
\begin{align*}
R(u) &= u= V(u^1)\over\cdots\over V(u^n), 
\end{align*}
and of course, being an algebra morphism, it sends the free products of 
$\Hr_Y$ into the over products of $\Ha$, that is 
\begin{align*}
R(u_1 \cdots u_m) &= u_1\over\cdots\over u_m. 
\end{align*}
To show the identity (\ref{Dr-Da}), it suffices to show that 
$(R\otimes R)\Dr_Y$ satisfies the same recursive relation (\ref{Da-def}) 
which defines $\Da$. 
For this purpose, we introduce a coaction $\dr_Y$ of $\Hr_Y$ on itself.  

Let us restrict the right action of $\Gd_Y(A)$ on $\Gr_Y(A)$ 
of Theorem~\ref{Gd-action-Glr} to the map 
$\Gr_Y(A)\times\Gr_Y(A)\longrightarrow\Gr_Y(A)$ given by 
$(\rho_f)^{\rho_g}=\rho_{f^{\rho_g}}$. We obtain a right action of $\Gr_Y(A)$ 
on itself. 
Its dual map on the coordinate rings can be found from the coaction 
$\ddif_Y$, given by Eq.~(\ref{ddif_P-def}) for $\P=\Dup$ and which 
in fact coincides with the coaction $\dinv_Y$, by applying the projection 
$P:\Hd_Y\longrightarrow\Hr_Y$ described in the proof of Theorem~\ref{Hr_P}. 
In conclusion, we obtain the map 
$\dr_Y:\Hr_Y\longrightarrow\Hr_Y\otimes\Hr_Y$ given on the generators 
$u\in\Y$ by 
\begin{align}
\dr_Y(u) &= (\Id\otimes P)\ddif_Y(u) 
= \sum_{u=\mu_t(s_1\over\treeAsmall,\dots,s_{|t|}\over\treeAsmall)}
t \otimes s_1\cdots s_{|t|}, 
\label{dr_Y(u)}
\end{align}
where, again, from now on we suppose that the sums run over all the trees 
in the set $\Y$ if they appear as subindeces of the tree-product $\mu$ 
(in this case $t$), and to the set $Y$ if they appear inside 
the arguments of $\mu$ or anywhere else (in this case $s_1,...,s_{|t|}$). 
Now we compute an explicit formula for $\Dr_Y(V(u))$ and for $\dr_Y(V(u))$, 
and show that relations (\ref{Da-def}) and (\ref{da-def}) are satisfied 
after projecting by $R$.

If $u=\treeO$, and $V(u)=\treeA$, we can easily compute 
\begin{align*}
\Dr_Y(\treeA)&=\treeA\otimes 1+1\otimes\treeA, \\ 
\dr_Y(\treeA)&=\treeA\otimes 1. 
\end{align*}
Therefore $\Dr_Y(\treeA)= 1\otimes\treeA+\dr_Y(\treeA)$. The 
relation (\ref{Da-def}) is satisfied in $\Hr_Y$, and therefore 
it is satisfyed after applying the algebra morphism $R$. 
Since $\dr_Y(\treeO)=\treeO\otimes\treeO$, and $\treeA=V(\treeO)$, 
relation (\ref{da-def}) is also satisfied in $\Hr_Y$. 

Now suppose that $u\neq\treeO$. 
In Eq.~(\ref{Dr_Y(u)}), we replace the tree $u$ by the tree $V(u)$, 
and obtain
\begin{align*}
\Dr_Y(V(u)) &= 1 \otimes V(u) + 
\sum_{V(u)=\mu_{\tilde{t}}
(\tilde{s}_1\over\treeAsmall,\dots,\tilde{s}_{|\tilde{t}|}\over\treeAsmall)
\over \tilde{s}_{|\tilde{t}|+1}}
\tilde{t} \otimes \tilde{s}_1\over\cdots\over\tilde{s}_{|\tilde{t}|} 
\over\tilde{s}_{|\tilde{t}|+1}.  
\end{align*}
Since $\tilde{t}\neq\treeO$, the tree
$\mu_{\tilde{t}}(\tilde{s}_1\over\treeA,\dots,\tilde{s}_{|\tilde{t}|}
\over\treeA)\over \tilde{s}_{|\tilde{t}|+1}$ can be of the form
$V(u)$ only if $\tilde{s}_{|\tilde{t}|+1}=\treeO$,
$\tilde{t}= V(t)$ with $t\in Y$, and $\tilde{s}_1=\treeO$.  
The case $t=\treeO$ corresponds to $u=\treeO$, and we exclude it.  
For $t\neq\treeO$, we write $V(t)=\mu_{\treeBA}(\treeA,t)$ and apply the 
associativity of $\mu$ to conclude that
\begin{align*}
V(u)&=\mu_{V(t)}
(\treeA,\tilde{s}_2\over\treeA,\dots,\tilde{s}_{|t|+1}\over\treeA)
\over\treeO \\
&= V(\mu_t(\tilde{s}_2\over\treeA,\dots,\tilde{s}_{|t|+1}\over\treeA)), 
\end{align*}
and therefore 
$u=\mu_t(\tilde{s}_2\over\treeA,\dots,\tilde{s}_{|t|+1}\over\treeA)$. 
By renaming the trees $\tilde{s}_i=s_{i-1}$, we finally obtain 
\begin{align}
\Dr_Y(V(u))  
&= 1 \otimes V(u) 
+ \sum_{u=\mu_t(s_1\over\treeAsmall,\dots,s_{|t|}\over\treeAsmall)}
V(t) \otimes s_1\cdots s_{|t|}. 
\label{Dr_Y(V(u))}
\end{align}
Similarly, if we replace the tree $u$ by the tree $V(u)$ in 
Eq.~(\ref{dr_Y(u)}), we obtain
\begin{align}
\dr_Y(V(u)) &= 
\sum_{u=\mu_t(s_1\over\treeAsmall,\dots,s_{|t|}\over\treeAsmall)}
V(t) \otimes s_1\cdots s_{|t|}. 
\label{dr_Y(V(u))}
\end{align}
Therefore we have 
$\Dr_Y(V(u))=1\otimes V(u)+\dr_Y(V(u))$, and consequently the 
relation (\ref{Da-def}) is fulfilled already in $\Hr_Y$. 

Using (\ref{dr_Y(u)}), Eq.~(\ref{dr_Y(V(u))}) can be written as 
\begin{align*}
\dr_Y(V(u)) &= (V\otimes\Id) 
\left(\sum_{u=\mu_t(s_1\over\treeAsmall,\dots,s_{|t|}\over\treeAsmall)}
t \otimes s_1\dots s_{|t|}\right) \\ 
&= (V\otimes\Id) \dr_Y(u). 
\end{align*}
Let us develop $\dr_Y(u)$, for $u=u^l\over V(u^r)$. 
\bigskip 

\noindent{\em Assume $u^r=\treeO$.\/} 
Let us start by considering the case $u=u^l\over\treeA$. 
We already computed $\dr_Y(\treeA)=\treeA\otimes\treeO$. 
The sum in (\ref{dr_Y(u)}) is over all trees $t\in \Y$ and 
$s_1,\dots,s_{|t|}\in Y$ such that 
$u=u^l\over\treeA=\mu_t(s_1\over\treeA,\dots,s_{|t|}\over\treeA)$.  
This equality is possible if and only if $t=t^l\over\treeA$ and 
$s_{|t|}=s_{|t|}^l\over\treeA$. 
Then we distinguish two possible cases, let us list the contributions 
to the whole sum coming from each of them. 
\bigskip 

\noindent{\em Case 1.\/} 
If $t^l\neq\treeO$, and therefore $|t|=|t^l|+1$, then 
$u^l=\mu_{t^l}(s_1\over\treeA,\dots,s_{|t^l|}\over\treeA)\over s_{|t^l|+1}^l$. 
If we rename $t^l=:w$, the contribution can be written as 
\begin{align*}
\sum_{u^l=\mu_{w}(s_1\over\treeAsmall,\dots,s_{|w|}\over\treeAsmall)}
w\over\treeA \otimes s_1\cdots s_{|w|}. 
\end{align*}

\noindent{\em Case 2.\/} 
If $t^l=\treeO$, then 
$u=u^l\over\treeA=\mu_{\treeAsmall}(s_1\over\treeA)=s_1\over\treeA$, 
and therefore $s_1=u^l$. We then have the contribution 
\begin{align*}
&\treeA\otimes u^l. 
\end{align*}
Summing up the two contributions, we obtain 
\begin{align*}
\dr_Y(u^l\over\treeA)    
&= \treeA \otimes u^l 
+ \sum_{u^l=\mu_{w}(s_1\over\treeAsmall,\dots,s_{|w|}\over\treeAsmall)}
w\over\treeA\otimes s_1\cdots s_{|w|}, 
\end{align*}
and therefore 
\begin{align*}
(R\otimes R)\dr_Y(u^l\over\treeA)    
&= \treeA \otimes u^l 
+ \sum_{u^l=\mu_{w}(s_1\over\treeAsmall,\dots,s_{|w|}\over\treeAsmall)}
w\over\treeA\otimes s_1\over\cdots\over s_{|w|} \\ 
&=\left(1\otimes u^l 
+ \sum_{u^l=\mu_{w}(s_1\over\treeAsmall,\dots,s_{|w|}\over\treeAsmall)}
w\otimes s_1\over\cdots\over s_{|w|}\right) \over 
\left(\treeA\otimes 1\right) \\ 
&=\left[(R\otimes R)\Dr_Y(u^l)\right]\over\left[(R\otimes R)\dr_Y(\treeA)\right]. 
\end{align*}
Therefore Equation~(\ref{da-def}) holds for $u=u^l\over V(\treeO)$. 
\bigskip 

\noindent{\em Assume $u^r \neq \treeO$.\/} 
The sum in (\ref{dr_Y(u)}) is over all trees $t\in \Y$ and 
$s_1,\dots,s_{|t|+1}\in Y$ such that 
$u=u^l\over V(u^r)=\mu_t(s_1\over\treeA,\dots,s_{|t|}\over\treeA)$.  
Since $t \neq\treeO$, we split $t=t^l\over V(t^r)$,
and distinguish three possible cases. 
Let us list the contributions to the sum coming from each of them.
\bigskip 

\noindent{\em Case 1.\/} 
If $t^l\neq\treeO$ and $t^r =\treeO$, then 
$t=t^l\over\treeA=\mu_{\treeAB}(t^l,\treeA)$ 
with $|t|=|t^l|+1$. Using the properties of the $\mu$ product, 
the equality becomes 
\begin{align*}
u^l\over V(u^r) 
&= \mu_{t^l\over\treeA}(s_1\over\treeA,\dots,s_{|t^l|+1}\over\treeA) 
= \mu_{t^l}(s_1\over\treeA,\dots,s_{|t^l|}\over\treeA) 
\over s_{|t^l|+1}\over\treeA \\
& = \mu_{t^l}(s_1\over\treeA,\dots,s_{|t^l|}\over\treeA) 
\over s_{|t^l|+1}\over V(\treeO), 
\end{align*}
which is impossible because $u^r\neq\treeO$. 
\bigskip 

\noindent{\em Case 2.\/} 
If $t^l=\treeO$ and $t^r\neq\treeO$, then $t=V(t^r)=\mu_{\treeBA}(\treeA,t^r)$ 
with $|t|=|t^r|+1$. 
Using the associativity of the $\mu$ product we get
\begin{align*}
\mu_{V(t^r)}(s_1\over\treeA,\dots,s_{|t^r|+1}\over\treeA) &= 
s_1\over\treeA\under \mu_{t^r}(s_2\over\treeA,\dots,s_{|t^r|+1}\over\treeA). 
\end{align*}
This tree can be equal to $u=u^l\over\treeA\under u^r$ if and only if $s_1=u^l$
and $\mu_{t^r}(s_2\over\treeA,\dots,s_{|t^r|+1}\over\treeA)=u^r$. 
Let us rename the free trees as follows:  
$t^r=:y$, $s_i=:z_{i-1}$ for $i=2,\dots,|t^r|+1$. Then the contribution 
to the whole sum for this case is 
\begin{align*}
& \sum_{u^r=\mu_y(z_1\over\treeAsmall,\dots,z_{|y|}\over\treeAsmall)}
V(y) \otimes u^l z_1\cdots z_{|y|}. 
\end{align*}

\noindent{\em Case 3.\/} 
If $t^l\neq\treeO$ and $t^r\neq\treeO$, using again the properties of the
$\mu$ product we have  
\begin{align*}
u^l\over V(u^r) 
&= \mu_{t^l}(s_1\over\treeA,\dots,s_{|t^l|}\over\treeA) \over 
\mu_{V(t^r)}(s_{|t^l|+1}\over\treeA,\dots,s_{|t^l|+|t^r|+1}\over\treeA) \\ 
&= \mu_{t^l}(s_1\over\treeA,\dots,s_{|t^l|}\over\treeA) \over s_{|t^l|+1}\over
V\big(\mu_{t^r}(s_{|t^l|+2}\over\treeA,\dots,s_{|t^l|+|t^r|+1}\over\treeA)\big).
\end{align*}
Renaming the free trees as $t^l=:w$, $s_i=:v_i$ for $i=1,\dots,|t^l|+1$,
and $t^r=:y$, $s_{|t^l|+j}=:z_{j-1}$ for $j=2,\dots,|t^r|+1$, we obtain 
the following contribution to the sum: 
\begin{align*}
\underset{u^r=\mu_y(z_1\over\treeAsmall,\dots,z_{|y|}\over\treeAsmall)}
{\sum_{u^l=\mu_w(v_1\over\treeAsmall,\dots,v_{|w|}\over\treeAsmall)
\over v_{|w|+1}}}
w \over V(y) \otimes v_1\cdots v_{|w|+1}\ z_1\cdots z_{|y|}. 
\end{align*}
Summing up the two contributions, we obtain 
\begin{align*}
\dr_Y(u^l\over V(u^r)) &= 
\sum_{u^r=\mu_y(z_1\over\treeAsmall,\dots,z_{|y|}\over\treeAsmall)}
V(y) \otimes u^l z_1\cdots z_{|y|} \\ 
& \qquad 
+ \underset{u^r=\mu_y(z_1\over\treeAsmall,\dots,z_{|y|}\over\treeAsmall)}
{\sum_{u^l=\mu_w(v_1\over\treeAsmall,\dots,v_{|w|}\over\treeAsmall)
\over v_{|w|+1}}}
w \over V(y) \otimes v_1\cdots v_{|w|+1}\ z_1\cdots z_{|y|}. 
\end{align*}
Therefore 
\begin{align*}
(R\otimes R)\dr_Y(u^l\over V(u^r)) &= 
\sum_{u^r=\mu_y(z_1\over\treeAsmall,\dots,z_{|y|}\over\treeAsmall)}
V(y) \otimes u^l \over z_1\over\cdots\over z_{|y|} \\ 
& + \underset{u^r=\mu_y(z_1\over\treeAsmall,\dots,z_{|y|}\over\treeAsmall)}
{\sum_{u^l=\mu_w(v_1\over\treeAsmall,\dots,v_{|w|}\over\treeAsmall)
\over v_{|w|+1}}}
w \over V(y) \otimes v_1\over\cdots\over v_{|w|+1} \over 
z_1\over\cdots\over z_{|y|} \\ 
&= \left[ 1\otimes u^l 
+ \sum_{u^l=\mu_w(v_1\over\treeAsmall,\dots,v_{|w|}\over\treeAsmall)
\over v_{|w|+1}}
w \otimes v_1\over\cdots\over v_{|w|+1}\right] \over \\ 
&\hspace{5cm}
\sum_{u^r=\mu_y(z_1\over\treeAsmall,\dots,z_{|y|}\over\treeAsmall)}
V(y) \otimes z_1\over\cdots\over z_{|y|} \\ 
&= \left[(R\otimes R)\Dr_Y(u^l)\right]\over
\left[(R\otimes R)\dr_Y(V(u^l))\right]. 
\end{align*}
Therefore the recursive relation~(\ref{da-def}) holds for any 
$u=u^l\over V(u^r)$. 
\end{proof}
\bigskip 

In \cite{BFqedtree}, it was proved that there exists a non-commutative lift 
of $\Ha$, let us denote it by $\Hanc$. 
As a corollary to the previous result, we obtain a non-recursive formula 
for the coproduct $\Da$, which is still valid on $\Hanc$. 

\begin{corollary}
The free associative algebra 
$\Hanc=\Q\langle V(u),\ u\in Y\rangle$ and its abelian quotient 
$\Ha=\Q[V(u),\ u\in Y]$ are graded and connected Hopf algebras with 
coproduct defined on the generators by 
\begin{align*}
\Da(V(u))  
&= 1 \otimes V(u) 
+ \underset{u=\mu_t(s_1\over\treeAsmall,\dots,s_{|t|}\over\treeAsmall)}
{\underset{s_1,\dots,s_{|t|}\in Y} {\sum_{t\in \Y}}} \ \  
V(t) \otimes s_1\over\cdots\over s_{|t|}. 
\end{align*}
\end{corollary}
\bigskip 

To conclude, we apply the results of Section~\ref{section-action} on the 
semi-direct coproduct Hopf algebras to the Hopf algebras $\Ha$ and $\He$. 
The ``QED renormalization Hopf algebra'' $\Ha\ltimes\He$, introduced in 
\cite{BFqedtree}, is then the non-commutative lift of the coordinate ring 
of the group functor $\Ga\ltimes\Ge$. 

Then, note that the order map $\pi$ gives a surjective group morphism from 
$\Ga(A)$ to $\Gd(A)$. In fact, since $\Ga(A)$ is a subgroup of $\Gr_Y(A)$, 
and $\pi$ is a group homomorphism from $\Gr_Y(A)$ to $\Gd(A)$, 
it suffices to verify that $\pi$ is still surjective when restricted to 
$\Ga(A)$. This follows from the fact that $\Ga(A)$ contains $\Gd(A)$ via 
the inclusion $\iu$ of Proposition~\ref{Gd_P->Gd}, which is a section of 
$\pi$. 

Therefore, for any fixed algebra $A$, the ``QED renormalization group'' 
$\Ga(A)\ltimes\Ge(A)$ is projected to the semi-direct product of 
usual series $\Gd(A)\ltimes\Gi(A)$.


\bibliographystyle{alpha}

\end{document}